\documentclass[12pt]{article}
\usepackage{amssymb, amsmath, amsthm, amscd}
\usepackage[pdftex]{graphics}
\usepackage[utf8]{inputenc}
\usepackage[all,cmtip]{xy}
\usepackage{bbm}
\usepackage{enumitem}
\usepackage{setspace}
\usepackage{mathdots}

\usepackage[mathscr]{euscript}
\usepackage[colorlinks=true]{hyperref}
\hypersetup{colorlinks   = true}
\hypersetup{linkcolor=blue}
\usepackage{tikz}
\usetikzlibrary{cd}

\begin{document}

\mathchardef\mhyphen="2D
\newtheorem{The}{Theorem}[section]
\newtheorem{Lem}[The]{Lemma}
\newtheorem{Prop}[The]{Proposition}
\newtheorem{Cor}[The]{Corollary}
\newtheorem{Rem}[The]{Remark}
\newtheorem{Obs}[The]{Observation}
\newtheorem{SConj}[The]{Standard Conjecture}
\newtheorem{Titre}[The]{\!\!\!\! }
\newtheorem{Conj}[The]{Conjecture}
\newtheorem{Question}[The]{Question}
\newtheorem{Prob}[The]{Problem}
\newtheorem{Def}[The]{Definition}
\newtheorem{Not}[The]{Notation}
\newtheorem{Claim}[The]{Claim}
\newtheorem{Conc}[The]{Conclusion}
\newtheorem{Ex}[The]{Example}
\newtheorem{Fact}[The]{Fact}
\newtheorem{Formula}[The]{Formula}
\newtheorem{Formulae}[The]{Formulae}
\newtheorem{The-Def}[The]{Theorem and Definition}
\newtheorem{Prop-Def}[The]{Proposition and Definition}
\newtheorem{Lem-Def}[The]{Lemma and Definition}
\newtheorem{Cor-Def}[The]{Corollary and Definition}
\newtheorem{Conc-Def}[The]{Conclusion and Definition}
\newtheorem{Terminology}[The]{Note on terminology}
\newcommand{\C}{\mathbb{C}}
\newcommand{\R}{\mathbb{R}}
\newcommand{\N}{\mathbb{N}}
\newcommand{\Z}{\mathbb{Z}}
\newcommand{\Q}{\mathbb{Q}}
\newcommand{\Proj}{\mathbb{P}}
\newcommand{\Rc}{\mathcal{R}}
\newcommand{\Oc}{\mathcal{O}}
\newcommand{\Vc}{\mathcal{V}}
\newcommand{\Id}{\operatorname{Id}}
\newcommand{\pr}{\operatorname{pr}}
\newcommand{\rk}{\operatorname{rk}}
\newcommand{\del}{\partial}
\newcommand{\delbar}{\bar{\partial}}
\newcommand{\Cdot}{{\raisebox{-0.7ex}[0pt][0pt]{\scalebox{2.0}{$\cdot$}}}}
\newcommand\nilm{\Gamma\backslash G}
\newcommand\frg{{\mathfrak g}}
\newcommand{\fg}{\mathfrak g}
\newcommand{\Oh}{\mathcal{O}}
\newcommand{\Kur}{\operatorname{Kur}}
\newcommand\gc{\frg_\mathbb{C}}
\newcommand\jonas[1]{{\textcolor{green}{#1}}}
\newcommand\luis[1]{{\textcolor{red}{#1}}}
\newcommand\dan[1]{{\textcolor{blue}{#1}}}

\begin{center}

{\Large\bf SKT Hyperbolic and Gauduchon Hyperbolic Compact Complex Manifolds}

\end{center}

\begin{center}
{\large Samir Marouani }

\end{center}

\vspace{1ex}

\noindent{\small{\bf Abstract.} We introduce two notions of hyperbolicity for not necessarily K\"ahler even balanced $n$-dimensional compact complex manifolds $X$. The first, called {\it SKT hyperbolicity}, generalises Gromov's K\"ahler hyperbolicity by means of SKT metrics.  The second, called {\it Gauduchon hyperbolicity} by means of Gauduchon metrics. Our first main result in this paper asserts that every SKT hyperbolic $X$ is also Kobayashi/Brody hyperbolic and every Gauduchon hyperbolic $X$ is divisorially hyperbolic. The second main result is to prove a vanishing theorem for the $L^2$ harmonic spaces on the universal cover of a SKT hyperbolic manifold.

\vspace{2ex}
\section{Introduction}\label{section:Introduction}
S. Kobayashi called a complex manifold $X$, that need not be either K\"ahler or compact, {\it hyperbolic} if the pseudo-distance he had introduced on $X$ is actually a distance. Using the mapping decreasing proprety of this distance, one can show that every holomorphic map from the complex plane $\C$ to a {\it Kobayashi hyperbolic} manifold is {\it constant}. Conversely, Brody observed that a compact complex manifold $X$ is Kobayashi hyperbolic if every holomorphic map from $\C$ to $X$ is constant. The long standing Kobayashi-Lang conjecture predicts that, for a compact K\"ahler manifold $X$, if $X$ is Kobayashi hyperbolic then its canonical bundle $K_X$ is {\it ample}.\\

M. Gromov introduced in one of his seminal papers [Gro91], the notion of {\it K\"ahler hyperbolicty} for a compact K\"ahler manifold $X$. The manifold $X$ is called {\it K\"ahler hyperbolic} if $X$ admits a K\"ahler metric $\omega$ whose lift $\widetilde\omega$ to the universal cover $\tilde{X}$ of $X$ can be expressed as $$\widetilde\omega = d\alpha$$ for a {\it bounded} $1$-form $\alpha$ on $\tilde{X}$. As pointed out by Gromov, it is not hard to see that the K\"ahler hyperbolicity implies the Kobayashi hyperbolicity.\\
The K\"ahler hyperbolicity is generalized in [MP22a] to what we call {\bf balanced hyperbolicity}. This is done by replacing the K\"ahler metric in the K\"ahler hyperbolicity by a {\it balanced metric}.
Meanwhile, a compact complex $n$-dimensional manifold $X$ is said to be balanced hyperbolic if it carries a balanced metric $\omega$  such that $\omega^{n-1}$ is $\widetilde{d}$-bounded.
The Brody hyperbolicity is replaced by what we call {\bf divisorial hyperbolicity}. A compact complex manifold $X$ is called {\it divisorially hyperbolic} if there exists no non-trivial holomorphic map from $\C^{n-1}$ to $X$ satisfying certain {\it subexponential volume growth condition}. Where the main result in  [MP22a] asserts that every balanced hyperbolic $X$ is also divisorially hyperbolic (Thm 2.8), and therfore the following implication holds:

\vspace{3ex}  

\hspace{12ex} $\begin{array}{lll} X \hspace{1ex} \mbox{is {\bf K\"ahler hyperbolic}} & \implies & X \hspace{1ex} \mbox{is {\bf Kobayashi/Brody hyperbolic}} \\
 \rotatebox{-90}{$\implies$} &  & \rotatebox{-90}{$\implies$} \\
 X \hspace{1ex} \mbox{is {\bf balanced hyperbolic}} & \implies & X \hspace{1ex} \mbox{is {\bf divisorially hyperbolic}} \\
 \rotatebox{90}{$\implies$} &  &   \\
 X \hspace{1ex} \mbox{is {\bf degenerate balanced}} & &\end{array}$\\
 
\vspace{3ex}  
In this paper we present a hyperbolicity theory where the K\"ahler metric is replaced by the SKT metric and the balanced metric is replaced by the Gauduchon metric on $n$-dimensional compact complex manifolds. The notions we introduce are weaker, thus more inclusive, than their classical counterparts. In particular, the setting does not need to be K\"ahler or even balanced. Our motivation stems from the existence of many interesting examples of non-K\"ahler compact complex manifolds that exhibit hyperbolicity features in a generalized sense, which we now set out to explain.
Our tow first  main result in  this paper asserts that every SKT hyperbolic $X$ is also Kobayashi hyperbolic, and every Gauduchon hyperbolic $X$ is Divisorially hyperbolic, and therfore the following implication holds:
 
\vspace{3ex} 
\begin{The}\label{The:introd_hyperbolicity-implications} Let $X$ be a compact complex manifold. The following implications hold:

\vspace{3ex}  

 $\begin{array}{lll} X \hspace{1ex} \mbox{is {\bf K\"ahler hyperbolic}} \implies  X \hspace{1ex} \mbox{is {\bf SKT hyperbolic}} & \implies & X \hspace{1ex} \mbox{is {\bf Kobayashi hyperbolic}} \\
 \rotatebox{-90}{$\implies$} &  & \rotatebox{-90}{$\implies$} \\
 X \hspace{1ex} \mbox{is {\bf balanced hyperbolic}} \implies  X \hspace{1ex} \mbox{is {\bf Gauduchon hyperbolic}}& \implies & X \hspace{1ex} \mbox{is {\bf divisorially hyperbolic}} \\
 \rotatebox{90}{$\implies$} &  &   \\
 X \hspace{1ex} \mbox{is {\bf degenerate balanced}} & &\end{array}$ 

\end{The}
\vspace{3ex} 
Our second main result is to prove a vanishing theorem for the $L^2$ harmonic spaces on the universal cover of a {\bf SKT hyperbolic} manifold. With conditions less stringent than K\"ahler hyperbolicity, we obtain the same result as in 1.4.A. Theorem in [Gro91].
\section{Aeppli Cohomology and SKT Metrics}

Given a compact complex $n$-dimensional manifold $X$, recall that the Bott-Chern and Aeppli cohomology groups of any bidegree $(p,\,q)$ of $X$ are classically defined, using the spaces $C^{r,\,s}(X) = C^{r,\,s}(X,\,\C)$ of smooth $\C$-valued $(r,\,s)$-forms on $X$, as \begin{eqnarray*}H^{p,\,q}_{BC}(X,\,\C) & = & \frac{\ker(\partial:C^{p,\,q}(X)\to C^{p+1,\,q}(X))\cap\ker(\bar\partial:C^{p,\,q}(X)\to C^{p,\,q+1}(X))}{\mbox{Im}\,(\partial\bar\partial:C^{p-1,\,q-1}(X)\to C^{p,\,q}(X))} \\
  H^{p,\,q}_A(X,\,\C) & = & \frac{\ker(\partial\bar\partial:C^{p,\,q}(X)\to C^{p+1,\,q+1}(X))}{\mbox{Im}\,(\partial:C^{p-1,\,q}(X)\to C^{p,\,q}(X)) + \mbox{Im}\,(\bar\partial:C^{p,\,q-1}(X)\to C^{p,\,q}(X))}.\end{eqnarray*}
The \mbox{Bott-Chern} Laplacian $\Delta_{BC}$ and the \mbox{Aeppli} Laplacian $\Delta_{A}$ are the $4^{th}$ order elliptic
differential operators defined respectively as
\begin{eqnarray}
    \Delta_{BC}:=\partial^\star\partial+\bar\partial^\star\bar\partial+(\partial\bar\partial)^\star(\partial\bar\partial)+(\partial\bar\partial)(\partial\bar\partial)^\star+(\partial^\star\bar\partial)^\star(\partial^\star\bar\partial)+(\partial^\star\bar\partial)(\partial^\star\bar\partial)^\star,
\end{eqnarray}
and
\begin{eqnarray}
    \Delta_A:=\partial\partial^\star+\bar\partial\bar\partial^\star+(\partial\bar\partial)^\star(\partial\bar\partial)+(\partial\bar\partial)(\partial\bar\partial)^\star+(\partial\bar\partial^\star)(\partial\bar\partial^\star)^\star+(\partial\bar\partial^\star)^\star(\partial\bar\partial^\star),
\end{eqnarray}
where $\partial^\star=-\star\bar\partial\star,$ and $\star=\star_\omega:\Lambda^{p,q}T^\star X\to\Lambda^{n-q,n-p}T^\star X$ is the Hodge-star isomorphism defined by
$\omega$ for arbitrary $p, q = 0, . . . , n$. Note that  $\star\Delta_{BC}=\Delta_A\star$ and $\Delta_{BC}\star=\star\Delta_A$, then $u \in \ker\Delta_{BC}$ if and only if $\star u\in \ker \Delta_A$.\\
The  Bott-Chern Laplacian is elliptic and formally self-adjoint, so it induces a three-space decomposition
\begin{eqnarray}
    C^\infty_{p,q}= \ker \Delta_{BC} \oplus \mbox{Im}\partial\bar\partial\oplus(\mbox{Im}\partial^\star+\mbox{ Im}\bar\partial^\star)
\end{eqnarray}
that is orthogonal w.r.t. the $L^2$ scalar product defined by $\omega$. We have
\begin{eqnarray}
\ker\partial\cap\ker\bar\partial=\ker\Delta_{BC}\oplus\mbox{Im}\partial\bar\partial
.\end{eqnarray}
We also have \begin{eqnarray}
    \mbox{Im}\Delta_{BC}=\mbox{Im}\partial\bar\partial\oplus(\mbox{Im}\partial^\star+\mbox{Im}\bar\partial^\star).
\end{eqnarray}
Similarly, the $4^{th}$ order Aeppli Laplacian is elliptic and formally self-adjoint, so it induces a three-space decomposition
\begin{eqnarray}
     C^\infty_{p,q}= \ker\Delta_{A} \oplus \mbox{Im}(\partial\bar\partial)^\star\oplus(\mbox{Im}\partial+\mbox{Im}\bar\partial),
\end{eqnarray}
that is orthogonal w.r.t. the $L2$
scalar product defined by $\omega$. We have
\begin{eqnarray}\label{Hodge_decomp_Delta_A}
    \ker(\partial\bar\partial)=\ker\Delta_A\oplus(\mbox{Im}\partial+\mbox{Im}\bar\partial).
\end{eqnarray}
We also have
\begin{eqnarray}
    \mbox{Im}\Delta_A=\mbox{Im}(\partial\bar\partial)^\star\oplus(\mbox{Im}\partial+\mbox{Im}\bar\partial).
\end{eqnarray}

\subsection{Cones of classes of metrics}
Recall a classical notion introduced by Popovici [Pop15]: the Gauduchon cone ${\cal G}_X$ of a compact complex manifold $X$ is the set of all
Bott-Chern cohomology classes of Gauduchon metrics on $X$. As such, it is an open convex cone in $H^{n-1, n-1}_{BC}(X, \R)$, that generalises the K\"ahler cone ${\cal K}_X$. Moreover, It is never
empty,  thanks to the existence of the Gauduchon metric on all
compact complex manifold. \\
Recall that a Hermitian metric $\omega$ on $X$ is said to be a {\it SKT metric} if $\partial\bar\partial\omega=0$. Any such metric $\omega$, defines an Aeppli cohomology class associated with $\omega$.
We start by introducing the following analogue of the Gauduchon cone in bidegree $(1,1)$.
\begin{Def}  Let $X$ be a compact complex manifold with $\dim_\C X =n$.
The {\bf SKT cone} of $X$ is the set : $${\cal SKT}_X:=\bigg\{[\omega]_A\in H^{1,\,1}_A(X,\,\R)\,\mid\,\omega\hspace{1ex}\mbox{is a SKT metric}\hspace{1ex} \mbox{on}\hspace{1ex}X\bigg\}\subset H^{1,\,1}_A(X,\,\R).$$
Any element $[\omega]_A$
 of the SKT cone ${\cal SKT}_X$ is called an {\bf Aeppli-SKT classe}.
\end{Def}
 We will also introduce  the {\it cone in cohomology of bidegree $(n-1, n-1)$} of $X$ which generalises the  the pseudo-effective cone ${\cal E}_X$ of Bott-Chern cohomology classes of $d$-closed semi-positive $(1, 1)$-currents introduced by Demailly [Dem92]. 
\begin{Def} Let $X$ be a compact complex manifold. The {\bf co-pseudo-effective cone}  of $X$ is the set
    
$${\cal E}^{n-1,n-1}_X:=\bigg\{[T]_{BC}\in H^{n-1,\,n-1}_{BC}(X,\,\R)\,\slash\, T\geq 0\hspace{1ex} d\mbox{-closed} \hspace{1ex} (n-1,\,n-1)\mbox{-current} \hspace{1ex} \mbox{on} \hspace{1ex} X\bigg\}.$$
Any element $[T]_{BC}$ of the co-pseudo-effective cone ${\cal E}^{n-1,n-1}_X$ is called a {\bf co-pseudo-effective class}.
\end{Def} 

\vspace{2ex}
 \begin{Prop}  Let $X$ be a compact complex manifold. Then
 \begin{enumerate}
     \item[(i)] The {\it SKT  cone } ${\cal SKT}_X$ is an open convex subset of $H^{1,1}_A(X,\,\R)$.
     \item[(ii)] The set ${\cal E}^{n-1,n-1}_X$ is a closed convex cone in $ H^{n-1,\,n-1}_{BC}(X,\,\R)$.
     \item[(iii)]   The following two statements hold
     \begin{enumerate}
         \item[(a)]
   Given any class $\mathfrak{c}^{n-1,\,n-1}_{BC}\in H^{n-1,\,n-1}_{BC}(X,\,\R)$, the following equivalence holds: $$\mathfrak{c}^{n-1,\,n-1}_{BC}\in{\cal E}^{n-1,n-1}_X\iff\mathfrak{c}^{n-1,\,n-1}_{BC}.\mathfrak{c}^{1,\,1}_A\geq 0 \hspace{2ex} \mbox{for every class} \hspace{1ex} \mathfrak{c}^{1,\,1}_A\in{\cal SKT}_X.$$

\item[(b)] Given any class $\mathfrak{c}^{1,\,1}_A\in H^{1,\,1}_A(X,\,\R)$, the following equivalence holds: $$\mathfrak{c}^{1,\,1}_A\in\overline{\cal SKT}_X\iff\mathfrak{c}^{n-1,\,n-1}_{BC}.\mathfrak{c}^{1,\,1}_A\geq 0 \hspace{2ex} \mbox{for every class} \hspace{1ex} \mathfrak{c}^{n-1,\,n-1}_{BC}\in{\cal E}^{n-1,n-1}_X.$$
  \end{enumerate}
 \end{enumerate}
 \end{Prop}
 \noindent {\it Proof.}
 
 (i)\ If $\omega_1$ and $\omega_2$ are {\bf SKT} metrics on $X$, so is any linear combination
$\lambda\omega_1+\mu\omega_2$ with $\lambda,\, \mu$  non-negative reals. Therefore, ${\cal SKT}_X$ is a convex cone. To show that ${\cal SKT}_X$ is an open subset of $H^{1,1}_A(X,\, \R)$, let us equip the finite-dimensional vector space $H
^{1, 1}_A(X,\, \R) $ with
an arbitrary norm $||\quad||$ (e.g. the Euclidian norm after we have fixed a basis ; at
any rate, all the norms are equivalent). Let $[\omega]_A\in{\cal SKT}_X$ be an arbitrary element, where $\omega > 0$ is some SKT metric on $X$. Let $\alpha\in H
^{1, 1}_A(X,\, \R) $
be a class such that $||\alpha - [\omega]_A|| < \epsilon$ for some small $ \epsilon > 0$. Fix any Hermitian metric $\widetilde\omega$ on $X$ and consider the Aeppli Laplacian $\Delta_A$ defined by
$\widetilde\omega$ inducing the Hodge isomorphism $H
^{1, 1}_A(X,\, \R) \simeq H^{1,1}_{\Delta_A}
(X,\R)$. Let
$\alpha_0 \in H
^{1, 1}_{\Delta_A}
(X, \R)$ be the $\Delta_A$-harmonic representative of the class $\alpha$.
 Since $\omega\in\mbox {ker}(\partial\bar\partial)$, (\ref{Hodge_decomp_Delta_A})  
 gives a unique decomposition
$$ \omega= \omega_0 + (\partial u + \bar\partial v) \mbox{ with } \Delta_A\omega_0 = 0.$$
If we set $ \gamma:= \alpha_0 + (\partial u + \bar\partial v)$ (with the same forms $u, v$ as for $\omega$), then $\partial\bar\partial\gamma = 0$, $\gamma$ represents
the Aeppli-SKT class $\alpha$ and we have $$ ||\gamma -\omega||_{C^0} = ||\alpha_0-\omega_0||_{C^0} \leq C ||\alpha -[\omega]_A|| < C\epsilon,
$$
for some constant $C > 0$ induced by the Hodge isomorphism. (We have chosen the $C^0$ norm on $H
^{1, 1}_{\Delta_A}
(X, \R)$ only for the sake of convenience.) Thus, if $\epsilon > 0$ is chosen sufficiently small, the $(1, 1)$-form $\gamma$ must be positive definite since $\omega$
is. Thus $\gamma$ is a {\bf SKT} metric and 
represents the original
Aeppli-SKT class $\alpha$, so $\alpha\in{\cal SKT}_X$.

(ii)\  Let $T_k$ be closed positive $(n-1,n-1)$-currents such that the classes $\{T_k\}$ converge to a
limit $\{\Theta\}$. Then $\int_
XT_k\wedge\omega$
converges to $\int_X\Theta\wedge\omega$
. In particular the sequence
$T_k$ is bounded in mass, and therefore weakly compact. If $T_{k_\nu} \to T$ is a weakly
convergent sub-sequence, the limit $T$ is a closed positive current, so $\{T \} = \{\Theta\}$ is co-pseudo-effective. Therefore ${\cal E}^{n-1,n-1}_X$ is closed.

 (iii)\ By the Serre-type duality (see e.g. [Pop15]): \begin{equation}\label{eqn:BC-A_duality}H^{n-1,\,n-1}_{BC}(X,\,\C)\times H^{1,\,1}_A(X,\,\C)\longrightarrow\C, \hspace{3ex} (\{u\}_{BC},\,\{v\}_A)\mapsto\{u\}_{BC}.\{v\}_A:=\int\limits_Xu\wedge v,\end{equation}
we can obtain (a) and (b) as a reformulation of a result in Lamari's work from [Lam99].
 \hfill $\Box$
 \subsection{Some properties of SKT manifolds}
 
On a complex manifold $X$ with $\mbox{dim}_\C X=n$, we will often use the following standard formula (cf. e.g. [Voi02, Proposition 6.29, p. 150]) for the Hodge star operator $\star = \star_\omega$ of any Hermitian metric $\omega$ applied to $\omega$-{\it primitive} forms $v$ of arbitrary bidegree $(p, \, q)$: \begin{eqnarray}\label{eqn:prim-form-star-formula-gen}\star\, v = (-1)^{k(k+1)/2}\, i^{p-q}\,\omega_{n-p-q}\wedge v, \hspace{2ex} \mbox{where}\,\, k:=p+q.\end{eqnarray} Recall that, for any $k=0,1, \dots , n$, a $k$-form $v$ is said to be $(\omega)$-{\it primitive} if $\omega_{n-k+1}\wedge v=0$ and that this condition is equivalent to $\Lambda_\omega v=0$, where $\Lambda_\omega$ is the adjoint of the operator $\omega\wedge\cdot$ (of multiplication by $\omega$) w.r.t. the pointwise inner product $\langle\,\,,\,\,\rangle_\omega$ defined by $\omega$.

\begin{Lem}\label{Lem:1-forms_Delta-harm} Let $\omega$ be a Hermitian metric on a complex manifold $X$ with $\mbox{dim}_\C X=n$. Fix a primitive form $\phi \in C^\infty_{p,q}(X,\,\C)$ with $p+q=n-1$.

  \vspace{3ex}

If $\omega$ is {\bf SKT} and $\partial\phi=\bar\partial\phi=0$, then 
\begin{eqnarray}\label{complete_case}\partial^\star(\omega\wedge \phi) = 0,\quad\quad\bar\partial^\star(\omega\wedge\phi)=0,\end{eqnarray}
In particular if $X$ is compact, then 
\begin{eqnarray*}\Delta_A(\omega\wedge \phi) = 0,\end{eqnarray*} where $\Delta_A$ is the Aeppli  Laplacian   induced by $\omega$.
\vspace{3ex}

\end{Lem} 
\noindent {\it Proof.} 
Let $\phi$ be a primitive form, the standard formula (\ref{eqn:prim-form-star-formula-gen}) we get: $\star\, \phi = (-1)^\frac{n(n-1)}{2}i^{p-q}\omega\wedge \phi$, hence $\star\,(\omega\wedge \phi) = -(-1)^\frac{n(n-1)}{2}i^{q-p}\phi$. Meanwhile, $d^\star = -\star d\star$, so applying $-\star d$ to the previous identity, we get \begin{eqnarray*}d^\star(\omega\wedge \phi)=(-1)^\frac{k(k+1)}{2}i^{q-p}(\star\partial \phi +\star\bar\partial \phi).\end{eqnarray*} 
By assumption, we have $\partial\phi=\bar\partial\phi=0$, then, on the one hand we have \begin{eqnarray*}d^\star(\omega\wedge \phi) &=& 0\\ &\Longleftrightarrow& \\
\partial^\star (\omega\wedge \phi)=0 &\mbox{ and } &\bar\partial^\star (\omega\wedge \phi)=0 \end{eqnarray*} 
and on the other hand and since $\omega$ is {\bf SKT}, we get

 \begin{eqnarray*}
    \partial\bar\partial(\omega\wedge \phi)=(\partial\bar\partial\omega)\wedge \phi=0.
\end{eqnarray*}

The proof of lemma is complete since on a compact complex manifold $X$, $$\alpha\in\mbox{ Ker }(\Delta_A)\Longleftrightarrow \partial^{\star} \alpha=0,\quad\bar\partial^{\star} \alpha=0,\quad\partial\bar\partial\alpha=0$$.
\hfill $\Box$
\begin{Cor}
Let $(X,\omega)$ be a {\bf SKT} $n$-dimensional compact complex manifold and $\phi$ be a $(n-1,0)$-form (respectively $(0,n-1)$-form) on $X$ such that $\partial\phi=0$ and $\bar\partial\phi=0$. Then,
\begin{eqnarray}\Delta_A(\omega\wedge \phi) = 0,\end{eqnarray} where $\Delta_A$ is the Aeppli  Laplacian   induced by $\omega$.
\end{Cor}
\begin{Lem-Def}\label{Lem-Def:primitive-classes_bideg11} Let $\omega$ be a {\bf SKT} metric on a compact complex manifold $X$ with $\mbox{dim}_\C X=n$. The linear map: \begin{equation*}\label{eqn:primitive-classes_bideg11}[\omega]_A\wedge\cdot : H^{n-1,\,n-1}_{BC}(X,\,\C)\longrightarrow H^{n,\,n}_A(X,\,\C)\simeq\C, \hspace{3ex} [\Gamma]_{BC}\longmapsto[\omega\wedge\Gamma]_A,\end{equation*} is {\bf well defined} and depends only on the cohomology class $[\omega]_A\in H^{1,\,1}_A(X,\,\C)$. Moreover any SKT metric gives rise to a non-zero class in Aeppli cohomology. 
We set: \begin{equation*}H^{n-1,\,n-1}_{BC}(X,\,\C)_{\omega\mhyphen prim}:=\ker\bigg([\omega]_A\wedge\cdot\bigg)\subset H^{n-1,\,n-1}_{BC}(X,\,\C)\end{equation*} and we call its elements {\bf ($\omega$)-primitive} Bott-Chern $(n-1,\,n-1)$-classes.

\end{Lem-Def}
\noindent {\it Proof.} The well-definedness follows at once from the identities: \begin{eqnarray*}\partial\bar\partial(\omega\wedge\Gamma) & = & \partial\bar\partial\omega\wedge\Gamma = 0,  \hspace{5ex} \Gamma\in C^\infty_{n-1,\,n-1}(X,\,\C)\cap\ker d,\\
  \omega_{n-1}\wedge\partial\bar\partial\Phi & = & \partial(\omega\wedge\bar\partial\Phi) + \bar\partial(\Phi\wedge\partial\omega)\in\mbox{Im}\,\partial + \mbox{Im}\,\bar\partial, \hspace{5ex} \Phi\in C^\infty_{n-2,\,n-2}(X,\,\C),\end{eqnarray*} where the latter takes into account the fact that $\partial\bar\partial\omega=0$.

That the map $[\omega]_A\wedge\cdot$ depends only on the Aeppli cohomology class $[\omega]_A$ follows from: \begin{eqnarray*}(\omega + \partial\bar\alpha + \bar\partial\alpha)\wedge\Gamma - \omega\wedge\Gamma = \partial(\bar\alpha\wedge\Gamma) + \bar\partial(\alpha\wedge\Gamma)\in\mbox{Im}\,\partial + \mbox{Im}\,\bar\partial,\hspace{5ex} \Gamma\in C^\infty_{n-1,\,n-1}(X,\,\C)\cap\ker d.\end{eqnarray*}
Now we proof that any SKT metric gives rise to a non-zero class in Aeppli cohomology. Suppose there is a $(0,1)$-form $\alpha$ and a $(1,0)$-form $\beta$ such that $\omega=\partial\alpha+\bar\partial\beta$. 
We get \begin{eqnarray*} \big(d(\alpha+\beta)\big)^n =(\omega+\bar\partial\alpha+\partial\beta)^n&=&\sum_{k=0}^n\sum_{j=0}^k C_n^kC_j^k(\bar\partial\alpha)^j\wedge(\partial\beta)^{k-j}\wedge\omega^{n-k}\\ &=& \sum_{\{k=2j\}}\mathcal{C}_j(\bar\partial\alpha)^j\wedge(\partial\beta)^{j}\wedge\omega^{n-2j}. \end{eqnarray*} To show that $\sum_{\{k=2j\}}\mathcal{C}_j(\bar\partial\alpha)^j\wedge(\partial\beta)^{j}\wedge\omega^{n-2j}>0$, it suffices to check that the real form $(\bar\partial\alpha)^j\wedge(\partial\beta)^{j}\wedge\omega^{n-2j}$ is
weakly (semi)-positive at every point of $X$. (Recall that $\bar\partial\alpha$ is the conjugate of $\partial\beta$.) To this end,
note that the $(2j, 2j)$-form $(\bar\partial\alpha)^j\wedge(\partial\beta)^{j}$
is weakly semi-positive as the wedge product of a $(2, 0)$-
form and its conjugate (see [Dem97, Chapter III, Example 1.2]). Therefore, the $(n, n)$-form is (semi)-positive since the product of a weakly (semi)-positive form and a strongly (semi)-positive form is weakly (semi)-positive and $\omega$ is strongly positive (see [Dem97,
Chapter III, Proposition 1.11]). (Recall that in bidegrees $(0, 1)$, $(1,1)$ , $(n-1,n-1)$ and $(n, n)$, the notions of
weak and strong positivity coincide).
By Stokes, we get $$0<\int_X\sum_{\{k=2j\}}\mathcal{C}_j(\bar\partial\alpha)^j\wedge(\partial\beta)^{j}\wedge\omega^{n-2j}=\int_X d\big((\alpha+\beta)\wedge(d(\alpha+\beta)^{n-1}\big)=0.$$ Consequently \begin{eqnarray}\label{non_0_A_classe}  [\omega]_A\neq 0 .\end{eqnarray}
\hfill $\Box$
\begin{Lem}  Let $\omega$ be a {\bf SKT} metric on a compact complex manifold $X$ with $\mbox{dim}_\C X=n$. Then,
    $H^{n-1,\,n-1}_{BC}(X,\,\C)_{\omega\mhyphen prim}$ is a complex hyperplane of $H^{n-1,\,n-1}_{BC}(X,\,\C)$ depending only on the cohomology class $[\omega]_A$
\end{Lem}
\noindent {\it Proof.} 
Suppose that $H^{n-1,\,n-1}_{BC}(X,\,\C)_{\omega\mhyphen prim}= H^{n-1,\,n-1}_{BC}(X,\,\C)$. This translates to \begin{equation}\label{eqn:primitive_Aeppli-exact_char_proof_1}\omega\wedge\Gamma\in\mbox{Im}\,\partial + \mbox{Im}\,\bar\partial, \hspace{6ex}\forall\,\Gamma\in C^\infty_{n-1,\,n-1}(X,\,\C)\cap\ker d.\end{equation}
Since $\omega$ is $(\partial\bar\partial)$-closed, it has a unique $L^2_\omega$-orthogonal decomposition: \begin{equation}\label{Aeppli_exact}\omega = (\omega)_h + (\partial\bar\alpha_\omega + \bar\partial\alpha_\omega),\end{equation} with a $(1,\,1)$-form $(\omega)_h\in\ker\Delta_{A,\,\omega}$ and a $(1,\,0)$-form $\alpha_\omega$. 

On the other hand, $\alpha$ is $d$-closed, so it has a unique $L^2_\omega$-orthogonal decomposition: \begin{equation}\Gamma = \Gamma_h + \partial\bar\partial\Phi,\end{equation} where $\Gamma_h$ is $\Delta_{BC,\,\omega}$-harmonic and $\Phi$ is a smooth function on $X$. 

Thus, for every $\Gamma\in C^\infty_{n-1,\,n-1}(X,\,\C)\cap\ker d$, we get: \begin{eqnarray*}\omega\wedge\Gamma & = & (\omega)_h\wedge\Gamma + \partial(\bar\alpha_\omega\wedge\Gamma) + \bar\partial(\alpha_\omega\wedge\Gamma) \\
  & = & (\omega)_h\wedge\Gamma_h + \partial\bigg((\omega)_h\wedge\bar\partial\Phi\bigg) + \bar\partial\bigg(\Phi\,\partial(\omega)_h\bigg) + \partial(\bar\alpha_\omega\wedge\Gamma) + \bar\partial(\alpha_\omega\wedge\Gamma),\end{eqnarray*} where for the last identity we used the fact that $\bar\partial\partial(\omega)_h=0$.

Thanks to assumption (\ref{eqn:primitive_Aeppli-exact_char_proof_1}), the last identity implies that \begin{equation}\label{eqn:primitive_Aeppli-exact_char_proof_2}(\omega)_h\wedge\Gamma_h\in\mbox{Im}\,\partial + \mbox{Im}\,\bar\partial, \hspace{6ex}\forall\,\Gamma_h\in C^\infty_{n-1,\,n-1}(X,\,\C)\cap\ker\Delta_{BC,\,\omega}.\end{equation}

Meanwhile, since $(\omega)_h$ is $\Delta_{A,\,\omega}$-harmonic (and real), $\star_\omega(\omega)_h$ is $\Delta_{BC,\,\omega}$-harmonic (and real). Hence, \begin{equation*}\mbox{Im}\,\partial + \mbox{Im}\,\bar\partial\ni(\omega)_h\wedge\star_\omega(\omega)_h = |(\omega)_h|^2_\omega\,dV_\omega\geq 0,\end{equation*} where the first relation follows from (\ref{eqn:primitive_Aeppli-exact_char_proof_2}) by choosing $\Gamma_h = \star_\omega(\omega)_h$. Consequently, from Stokes's Theorem we get: $$\int\limits_X |(\omega)_h|^2_\omega\,dV_\omega = 0,$$ hence $(\omega)_h=0$, which implies by the $L^2$-orthogonal decomposition \ref{Aeppli_exact}, that $[\omega]_A=0$ which is absurd according (\ref{non_0_A_classe}).

 \hfill $\Box$

 \begin{Lem}\label{Lem:orthogonality_prim-class-rep_Gauduchon} Suppose there exists a {\bf SKT} metric $\omega$ on a compact complex manifold $X$. Then, for every $\Gamma\in C^\infty_{n-1,\,n-1}(X,\,\C)$ such that $d\Gamma=0$ and $[\Gamma]_{BC}\in H^{n-1,\,n-1}_{BC}(X,\,\C)_{\omega\mhyphen prim}$, we have: $$\langle\langle(\omega_{n-1})_h,\,\Gamma\rangle\rangle_\omega = 0,$$ where $\langle\langle\,\,,\,\,\rangle\rangle_\omega$ is the $L^2$ inner product induced by $\omega$.

\end{Lem}

\noindent {\it Proof.} Since $[\Gamma]_{BC}\in H^{n-1,n-1}_{BC}(X,\,\C)_{\omega\mhyphen prim}$ then, $\omega\wedge\Gamma$ is of the
shape $\partial \Phi+\bar\partial \Psi$ for some $\Phi\in C^\infty_{n-1,n}$ and $\Psi\in C^\infty_{n,n-1}$. By Stokes, we get: \begin{eqnarray*}\langle\langle\Gamma,\,(\omega_{n-1})_h\rangle\rangle_\omega & = & \int\limits_X\Gamma\wedge\star_\omega(\omega_{n-1})_h = \int\limits_X\Gamma\wedge(\omega)_h = \int\limits_X\Gamma\wedge(\omega  -\partial u-\bar\partial v) \\
  & = & \int\limits_X\Gamma\wedge\omega+\int\limits_X\partial(\Gamma\wedge u)+\int\limits_X\bar\partial(\Gamma\wedge v)= \int\limits_X\partial \Phi+\bar\partial \psi = 0,\end{eqnarray*}
  since $\partial\Gamma=0$, $\bar\partial\Gamma= 0$, $\partial(\Gamma\wedge u)=d(\Gamma\wedge u)$ and $\partial(\Gamma\wedge u)=d(\Gamma\wedge u)$.
  \hfill $\Box$
  
   \vspace{2ex}
   
  We shall now get a {\it Lefschetz-type decomposition} of $H^{n-1,n-1}_{BC}(X,\,\C)$, induced by an arbitrary SKT metric $\omega$, with $H^{n-1,n-1}_{BC}(X,\,\C)_{\omega\mhyphen prim}$ as a direct factor.
  \begin{Prop}
  Let $X$ be a {\bf SKT} compact complex manifold with $\mbox{dim}_\C X=n$. Then, the Bott-Chern cohomology space of bidegree $(n-1,n-1)$ has a {\bf Lefschetz-type} $L^2_\omega$-orthogonal {\bf decomposition}: \begin{equation}\label{eqn:Lefschetz-type-decomp_skt}H^{n-1,n-1}_{BC}(X,\,\C) = H^{n-1,n-1}_{BC}(X,\,\C)_{\omega\mhyphen prim}\oplus\C\cdot[(\omega_{n-1})_h]_{BC},\end{equation} where the $(\omega)$-primitive subspace $H^{n-1,n-1}_{BC}(X,\,\C)_{\omega\mhyphen prim}$ is a complex hyperplane of $H^{n-1,n-1}_{BC}(X,\,\C)$ depending only on the cohomology class $[\omega]_{A}\in H^{1,1}_{A}(X,\,\C)$, while $(\omega_{n-1})_h$ is the $\Delta_{BC,\omega}$-harmonic component of $\omega_{n-1}$ and the complex line $\C\cdot[(\omega_{n-1})_h]_{A}$ depends on the choice of the SKT metric $\omega$.
  \end{Prop}

  \begin{Lem}\label{Lem:Lefschetz-type-decomp_Gauduchon_class_coeff}  For every $\Gamma\in C^\infty_{n-1,\,n-1}(X,\,\C)\cap\ker d$, the coefficient of $[(\omega_{n-1})_h]_{BC}$ in the Lefschetz-type decomposition of $[\Gamma]_{BC}\in H^{n-1,\,n-1}_{BC}(X,\,\C)$ according to (\ref{eqn:Lefschetz-type-decomp_skt}), namely in \begin{equation}\label{eqn:Lefschetz-type-decomp_Aeppli-Gauduchon_class}[\Gamma]_{BC} = [\Gamma]_{BC,\,prim} + \lambda\,[(\omega_{n-1})_h]_{BC},\end{equation} is given by \begin{equation}\label{eqn:Lefschetz-type-decomp_Aeppli-Gauduchon_class_coeff}\lambda = \lambda_\omega([\Gamma]_{BC}) = \frac{[\omega]_A.[\Gamma]_{BC}}{||\omega_h||^2_\omega} = \frac{1}{||\omega_h||^2_\omega}\,\int\limits_X\Gamma\wedge\omega.\end{equation}

\end{Lem}

\noindent {\it Proof.} Since $[\Gamma]_{BC,\,prim}\in H^{n-1,n-1}_{BC}(X,\,\C)_{\omega\mhyphen prim}$, we have $[\omega]_A.[\Gamma]_{BC,\,prim}=0$, so \begin{eqnarray*}[\omega]_{A}.[\Gamma]_{BC} & = & \lambda\,\int\limits_X\omega\wedge(\omega_{n-1})_h = \lambda\,\int\limits_X\omega_h\wedge(\omega_{n-1})_h = \lambda\,\int\limits_X\omega_h\wedge\star_\omega(\omega)_h \\
  & = &  \lambda\,||(\omega)_h||^2_\omega .\end{eqnarray*} This gives (\ref{Lem:Lefschetz-type-decomp_Gauduchon_class_coeff}).  \hfill $\Box$
  \vspace{2ex}
  
Formula (\ref{eqn:Lefschetz-type-decomp_Aeppli-Gauduchon_class_coeff}) implies that $\lambda_\omega([\Gamma]_{BC})$ is {\it real} if the class $[\Gamma]_{BC}\in H^{n-1,\,n-1}_{BC}(X,\,\C)$ is real. Thus, we can define a {\it positive side} and a {\it negative side} of the hyperplane $$H^{n-1,\,n-1}_{BC}(X,\,\R)_{\omega\mhyphen prim}:=H^{n-1,\,n-1}_{BC}(X,\,\C)_{\omega\mhyphen prim}\cap H^{n-1,\,n-1}_{BC}(X,\,\R) \mbox{ in } H^{1,\,1}_{BC}(X,\,\R)$$ by \begin{eqnarray*}\label{eqn:positive-negative-sides_BC}\nonumber H^{n-1,\,n-1}_{BC}(X,\,\R)_\omega^{+} & := & \bigg\{[\Gamma]_{BC}\in H^{n-1,\,n-1}_{BC}(X,\,\R)\,\mid\,\lambda_\omega([\Gamma]_{BC})>0\bigg\}, \\
H^{n-1,\,n-1}_{BC}(X,\,\R)_\omega^{-} & := & \bigg\{[\Gamma]_{BC}\in H^{n-1,\,n-1}_{BC}(X,\,\R)\,\mid\,\lambda_\omega([\Gamma]_{BC})<0\bigg\}.\end{eqnarray*} These are open subsets of $H^{n-1,\,n-1}_{BC}(X,\,\R)$ that depend only on the cohomology class $[\omega]_A\in H^{1,\,1}_A(X,\,\R)$.

Since $[\Gamma]_{BC}$ is $\omega$-primitive if and only if $\lambda_\omega([\Gamma]_{BC})=0$, we get a {\it partition} of $H^{n-1,\,n-1}_{BC}(X,\,\R)$: \begin{eqnarray*}\label{eq:partition_H2_BC} H^{n-1,\,n-1}_{BC}(X,\,\R) = H^{n-1,\,n-1}_{BC}(X,\,\R)_\omega^{+}\cup H^{n-1,\,n-1}_{BC}(X,\,\R)_{\omega\mhyphen prim}\cup H^{n-1,\,n-1}_{BC}(X,\,\R)_\omega^{-}\end{eqnarray*} depending only on the cohomology class $[\omega]_A\in H^{1,\,1}_A(X,\,\R)$.

\vspace{2ex}

As a consequence of these considerations, we get
\begin{Prop}\label{Prop:psef-cone_positive-side-char} Let $X$ be a compact complex manifold with $\mbox{dim}_\C X=n$.\\ The {\bf  cone}  ${\cal E}^{n-1,n-1}_X\subset H^{n-1,\,n-1}_{BC}(X,\,\R)$ of $X$ is the intersection of all the {\bf non-negative sides} $$H^{n-1,\,n-1}_{BC}(X,\,\R)_\omega^{\geq 0}:=H^{n-1,\,n-1}_{BC}(X,\,\R)_\omega^{+}\cup H^{n-1,\,n-1}_{BC}(X,\,\R)_{\omega\mhyphen prim}$$  of hyperplanes $H^{n-1,\,n-1}_{BC}(X,\,\R)_{\omega\mhyphen prim}$ determined by cohomology classes $[\omega]_A\in H^{1,1}_A(X,\,\R)$: \begin{equation}\label{eqn:psef-cone_positive-side-char}{\cal E}^{n-1,n-1}_X = \bigcap\limits_{[\omega]_A\in H^{1,1}_A(X,\,\R) }H^{n-1,\,n-1}_{BC}(X,\,\R)_\omega^{\geq 0},\end{equation}  
\end{Prop}  
\noindent {\it Proof.}
 By the duality between the co-pseudo-effective cone ${\cal E}^{n-1,n-1}_X$ and the closure $\overline{\cal SKT}_X$ of the SKT cone, we know that a given class $[T]_{BC}\in H^{n-1,\,n-1}_{BC}(X,\,\R)$ lies in ${\cal E}^{n-1,n-1}_X$ (i.e. $[T]_{BC}$ can be represented by a closed {\it semi-positive} $(n-1,\,n-1)$-current) if and only if $$\int\limits_X T\wedge\omega\geq 0 \hspace{5ex} \mbox{for all}\hspace{1ex} [\omega]_A\in{\cal SKT}_X.$$ The last condition is equivalent to $\lambda_\omega([T]_{BC})\geq 0$, hence to $[T]_{BC}\in H^{n-1,\,n-1}_{BC}(X,\,\R)_\omega^{\geq 0}$, for all $[\omega]_A\in{\cal SKT}_X$, so the contention follows.
 \hfill $\Box$

\section{SKT and Gauduchon hyperbolicity}\label{section:bal-div_hyperbolicity}

In this section, we introduce and discuss two notions of hyperbolicity for not necessarily K\"ahler even balanced $n$-dimensional
compact complex manifolds that generalise Gromov's K\"ahler hyperbolicity and balanced hyperbolicity. 

\vspace{1ex}

Let $X$ be a compact complex manifold with $\mbox{dim}_\C X=n$. Fix an arbitrary Hermitian metric (i.e. a $C^\infty$ positive definite $(1,\,1)$-form) $\omega$ on $X$. The metric $\omega$ is said to be {\bf K\"ahler} $d\omega=0$. The manifold $X$ is said to be {\it K\"ahler} if it carries a {\it K\"ahler metric}. there are several well studied generalizations of the K\"ahlerness condition:
\vspace{1ex}
\begin{enumerate}

 \item[i)] $\omega $ is {\it balanced}, if $d\omega^{n-1}=0$.
\item[ii)] $\omega$ is {\it Gauduchon}, if $\bar\partial\partial\omega^{n-1}=0$, such a metric always exists on a compact complex manifold.
\item[iii)] $\omega$ is {\it strongly Gauduchon} if there is a $(n,n-2)$-form $\Gamma$ such that $\partial\omega^{n-1}=\bar\partial\Gamma$.
 \item[iv)] $\omega$ is Hermitian symplectic, if there is a $(0,2)$- form $\alpha$ such that $\partial\alpha=0$ and $\partial\omega=-\bar\partial\alpha$. 

  \item[v)] $\omega$ is {\bf SKT} (or pluriclosed), if $\partial\bar\partial\omega=0$
\end{enumerate}
Throughout the text, $\pi_X:\widetilde{X}\longrightarrow X$ will stand for the universal cover of $X$ and $\widetilde\omega=\pi_X^\star\omega$ will be the Hermitian metric on $\widetilde{X}$ that is the lift of $\omega$. Recall that a $C^\infty$ $k$-form $\alpha$ on $X$ is said to be $\widetilde{d}(\mbox{bounded})$ with respect to $\omega$ if $\pi_X^\star\alpha = d\beta$ on $\widetilde{X}$ for some $C^\infty$ $(k-1)$-form $\beta$ on $\widetilde{X}$ that is bounded w.r.t. $\widetilde\omega$. (See [Gro91] and [MP22]). In general, we propose the following definition which generalizes that of $\widetilde{d}$-bounded of a differential form.
\begin{Def}
a $C^\infty$ $k$-form $\phi$ on $X$ is said to be $\widetilde{(\partial+\bar\partial)}$-bounded with respect to $\omega$ if $\pi_X^\star\phi = \partial\alpha+\bar\partial\beta$ on $\widetilde{X}$ for some $C^\infty$ $(k-1)$-forms $\alpha$ and $\beta$ on $\widetilde{X}$ that are bounded w.r.t. $\widetilde\omega$.
    
\end{Def}

\vspace{1ex}

\subsection{ SKT hyperbolic manifolds}\label{subsection:bal-hyp_manifolds}
\hspace{1ex}
The first notion that we introduce in this work is the following
\hspace{1ex}
\begin{Def}\label{Def:bal-hyperbolic} Let $X$ be a compact complex manifold with $\mbox{dim}_\C X=n$. A Hermitian metric $\omega$ on $X$  is said to be {\bf SKT hyperbolic} if $\omega$ is SKT and  $(\widetilde{\partial+\bar\partial})-\mbox{bounded}$ with respect to $\omega$.

   The manifold $X$ is said to be SKT hyperbolic if it carries a {\it SKT hyperbolic metric}.

\end{Def}

Let us first notice the following implication: $$X \hspace{1ex} \mbox{is {\it K\"ahler hyperbolic}} \hspace{1ex} \implies \hspace{1ex} X \hspace{1ex} \mbox{is {\it SKT hyperbolic}}.$$

To see this, note the obvious fact that every K\"ahler metric is SKT and every $\widetilde{d}$-bounded form is in particular $(\widetilde{\partial+\bar\partial})-\mbox{bounded}$.

\begin{Lem}\label{Lem:powers_d-tilde-bounded} Let $(X,\,\omega)$ be a compact complex Hermitian manifold with $\mbox{dim}_\C X=n$. Let $k\in\{1,\dots , 2n\}$ and $\alpha\in C^\infty_k(X,\,\C)\cap\ker\partial\cap\ker\bar\partial$. If $\alpha$ is $(\widetilde{\partial+\bar\partial})$-bounded (with respect to $\omega$), then $\alpha^p$ is $(\widetilde{\partial+\bar\partial})$-bounded (with respect to $\omega$) for every non-negative integer $p$.

\end{Lem} 
\noindent {\it Proof.} By the $\widetilde{(\partial+\bar\partial)}-\mbox{boundedness}$ assumption on $\alpha$, $\pi_X^\star\alpha = \partial\beta+\bar\partial\gamma$ on $\widetilde{X}$ for some smooth $\widetilde\omega$-bounded $(k-1)$-forms $\beta$ and $\gamma$ on $\widetilde{X}$. Note that $\partial\beta+\bar\partial\gamma$ is trivially $\widetilde\omega$-bounded on $\widetilde{X}$ since it equals $\pi_X^\star\alpha$ and $\alpha$ is $\omega$-bounded on $X$ thanks to $X$ being {\it compact}. We get: \begin{equation*}\pi_X^\star\alpha^p = (\partial\beta+\bar\partial\gamma)\wedge(\pi_X^\star\alpha)^{p-1}= \partial(\beta\wedge(\pi_X^\star\alpha)^{p-1})+\bar\partial(\gamma\wedge(\pi_X^\star\alpha)^{p-1}), \end{equation*} where both $\beta$, $\gamma$ and $\pi_X^\star\alpha$ are $\widetilde\omega$-bounded, hence so is $(\partial\beta+\bar\partial\gamma)\wedge(\pi_X^\star\alpha)^{p-1}$.  \hfill $\Box$
\begin{Cor}
Let $(X,\,\omega)$ be a compact complex K\"aher manifold with $\mbox{dim}_\C X=n$. The following implication hold
$$ \omega \mbox{ is {\bf SKT} hyperbolic } \implies \omega^k \mbox{ is } (\widetilde{\partial+\bar\partial})\mbox{-bounded for all } k\in\{1,...,n\}$$

\end{Cor}
\vspace{3ex}
 
\begin{The}\label{The:bal-div-hyperbolic_implication} Every {\bf SKT hyperbolic} compact complex manifold is {\bf Kobayashi hyperbolic}.

\end{The}
  \noindent {\it Proof.} Let $X$ be a compact complex manifold, with $\mbox{dim}_\C X=n$, equipped with a {\it SKT hyperbolic} metric $\omega$. This means that, if $\pi_X:\widetilde{X}\longrightarrow X$ is the universal cover of $X$, we have $$\pi_X^\star\omega = \partial\alpha+\bar\partial\beta  \hspace{3ex} \mbox{on}\hspace{1ex} \widetilde{X},$$ where $\alpha$ and $\beta$ are a smooth, $\widetilde\omega$-bounded $(0,1)$-form and $(1,0)$-form  on $\widetilde{X}$ and $\widetilde\omega=\pi_X^\star\omega$ is the lift of the metric $\omega$ to $\widetilde{X}$.

Suppose there exists a  non-constant holomorphic map $f:\C\longrightarrow X$. We will prove that $f^\star\omega=0$ on $\C$, which implies  $f$ is a constant, a contradiction.

There exists a lift $\widetilde{f}$ of $f$ to $\widetilde{X}$, namely a holomorphic map $\tilde{f}:\C\longrightarrow\widetilde{X}$ such that $f=\pi_X\circ\tilde{f}$. The $(1,\,1)$-form $f^\star\omega$ is $\geq 0$ on $\C$ and may not be $>0$ on countably many points  $\Sigma\subset\C$. Thus, it can be written in the form (a) below for some $C^\infty$ function $\mu:\C\longrightarrow[0,\,+\infty)$: $$f^\star\omega\stackrel{(a)}{=} \mu(z)\,\frac{i}{2}\,dz\wedge d\bar{z}\stackrel{(b)}{=} \tilde{f}^\star(\pi_X^\star\omega)=\tilde{f}^\star(\partial\alpha+\bar\partial\beta)\stackrel{(c)}{=} d(\tilde{f}^\star(\alpha+\beta))  \hspace{3ex} \mbox{on}\hspace{1ex} \C,$$ where (b) follows from $f=\pi_X\circ\tilde{f}$, and (c) follows since $f^\star$ commutes with $d$
and by the fact that a
$(1,1)$-form is of maximal degree on $\C$.
  We have the following

  \begin{Claim}\label{Claim:f-tilde-star_gamma_bounded} The $1$-form $\tilde{f}^\star(\alpha+\beta)$ is $(f^\star\omega)$-bounded on $\C$.

  \end{Claim}

  \noindent {\it Proof of Claim.} For any tangent vector $v$, we have: \begin{eqnarray*}\big|\big(\tilde{f}^\star(\alpha+\beta)\big)(v)\big|^2 & = & \big|\big(\alpha+\beta\big)(\tilde{f}_\star v)\big|^2 \stackrel{(a)}{\leq} (C_1+C_2)\,|\tilde{f}_\star v|^2_{\widetilde\omega} \\
    & = & C\,|v|^2_{\tilde{f}^\star\widetilde\omega} \stackrel{(b)}{=} C\,|v|^2_{f^\star\omega},\end{eqnarray*} where $C>0$ is a constant independent of  $v$ that exists such that inequality (a) holds thanks to the {\it $\widetilde\omega$-boundedness} of $\alpha$ and $\beta$ on $\widetilde{X}$, while (b) follows from $\tilde{f}^\star\widetilde\omega = f^\star\omega$.   \hfill $\Box$

  \vspace{2ex}

  \noindent {\it End of Proof of Theorem \ref{The:bal-div-hyperbolic_implication}.} For any bounded open subset $\Omega\subset\C$, we let $\mbox{A}_\mu(\Omega)$ and $L_\mu(\partial\Omega)$ stand for the area  of $\Omega$ , respectively the length of $\partial\Omega$ w.r.t. the measure $f^\star\omega = \mu(z)\,(i/2)\,dz\wedge d\bar{z}$ on $\C$. For any $r>0$, letting $D_r:=\{z\in\C\,\mid\,|z|<r\}$ and $S_r:=\{z\in\C\,\mid\,|z|=r\}$, we get:

  \begin{eqnarray*}\mbox{A}_\mu(D_r) = \frac{i}{2}\int\limits_{D_r}\mu(z)\,dz\wedge d\bar{z} = \int\limits_0^r\bigg(\int\limits_{S_t}\mu\,d\sigma\bigg)\,dt  \stackrel{(a)}{\geq} \int\limits_0^r\bigg(\int\limits_{S_t}\mu^{\frac{1}{2}}\,d\sigma\bigg)^2\,\frac{1}{\bigg(\int\limits_{S_t}d\sigma\bigg)}\,dt,\end{eqnarray*} where $d\sigma$ is the length element of the circle  and inequality (a) follows from H\"older's inequality applied to the functions $\mu^{1/2}$ and $1$.

   Therefore, we get: \begin{eqnarray}\label{eqn:bal-div-hyperbolic_implication_proof_1_1}\mbox{A}_\mu(D_r) \geq \int\limits_0^r L_\mu^2\,(S_t)\,\frac{1}{2\pi\,t}\,dt.\end{eqnarray}

  Meanwhile, for every $r>0$, we also have: \begin{eqnarray}\label{eqn:bal-div-hyperbolic_implication_proof_2_1}\mbox{A}_\mu(D_r) = \int\limits_{D_r} f^\star\omega = \int\limits_{D_r}d(\tilde{f}^\star(\alpha+\beta)) = \int\limits_{S_r}\tilde{f}^\star(\alpha+\beta) \stackrel{(a)}{\leq} C\,\int\limits_{S_r}\mu^{\frac{1}{2}}\,d\sigma\stackrel{(b)}{=}C\,L_\mu(S_r),\end{eqnarray} where $C>0$ is a constant that exists such that inequality (a) holds thanks to Claim \ref{Claim:f-tilde-star_gamma_bounded} and equality (b) holds since the (1,1)-form $f^\star\omega$ already has maximal degree on $\C$ and is, therefore, proportional to $\frac{i}{2}dz\wedge d\bar{z}$ .

  Putting (\ref{eqn:bal-div-hyperbolic_implication_proof_1_1}) and (\ref{eqn:bal-div-hyperbolic_implication_proof_2_1}) together, we get: \begin{eqnarray}\label{eqn:bal-div-hyperbolic_implication_proof_3_1}\mbox{A}_\mu(D_r) \geq \frac{1}{2\pi\,C^{2}}\,\int\limits_0^r\frac{\mbox{A}_\mu(D_t)^2}{t}\,dt:=\frac{1}{2\pi\,C^2}\,F(r),\end{eqnarray} where the last equality is the definition of a function $F:(0,\,+\infty)\longrightarrow(0,\,+\infty)$.

  Deriving $F$, we get for every $r>0$: \begin{eqnarray*}t\,\frac{d}{dt}F(t) = \mbox{A}_\mu(D_t)^2 \geq 2\pi\,\frac{1}{4\pi^2C^2}(F(t))^2,\end{eqnarray*} where the last inequality follows from (\ref{eqn:bal-div-hyperbolic_implication_proof_3_1}). This amounts to \begin{eqnarray*}\frac{d}{dt}\bigg(\frac{-1}{F(t)}\bigg) \geq\frac{1}{4\pi^2C^2} .\frac{1}{t}, \hspace{3ex} t>0.\end{eqnarray*} Integrating the last identity over $t\in[a,\,b]$, with $0<a<b$ arbitrary, we get: \begin{eqnarray*}\frac{1}{F(a)}\geq\frac{1}{F(b)} + \frac{1}{4\pi^2C^2}\,\log\frac{b}{a}\geq\frac{1}{4\pi^2C^2}\,(\log b - \log a), \hspace{3ex} 0<a<b.\end{eqnarray*}

  For every fixed $a>0$, letting $b\to +\infty$ in the last identity, we get $F(a)=0$. This can only happen if $f^\star\omega=0$ on $\C$, which contradicts the fact that $f$ is non-constant on $\C$. \hfill $\Box$
  \begin{Prop}
Let $G$ be a complex Lie groups and $\Gamma$ a co-compact subgroup of $G$, then   there is no {\bf SKT hyperbolic metric} on $G/\Gamma$.
Moreover, $G/\Gamma$ can not be Kobayashi hyperbolic. 
\end{Prop}

  \noindent {\it Proof.}\footnote{The author is very grateful to M. Zaidenberg for 
 his collaboration to obtain this observation}

A Lie group G has 1-dimensional Lie subgroups $H$ passing to the neutral element $e$ in any tangent direction. The latter subgroup has for its universal covering a connected,  simply connected group of dimension $1$. Now, such a group is the complex line $\mathbb{ C}$. So, the universal covering is $\mathbb{ C} \rightarrow H$, and the $H$'s cover a  neighborhood $U$ of $e$. Kobayashi's metric vanishes along an $H$. Hence it is identically zero in $U$ because the distance from any point to $e$ is zero. In particular, the infinitesimal Kobayashi-Royden pseudometric on the tangent bundle $TU$ vanishes identically. We can cover $G$ by translates of $U$. Since the left multiplication is an isometry with respect to the Kobayashi and the Kobayashi-Royden pseudometics, the Kobayashi-Royden pseudometric of $G$ is identically zero. Then also the Kobayashi pseudometric is.
It follows that it is identically zero on any quotient $G / \Gamma$ where $\Gamma$ is a subgroup of $G$, that is, on any homogeneous space of $G$. 
\hfill $\Box$

\vspace{3ex}
We now adapt in a straightforward way to our context the first part of the proof of [CY18,
Proposition 2.11], where the non-existence of rational curves in compact K\"ahler hyperbolic manifolds
was proved, and get the same result with a weaker condition on $\omega$.
\begin{Prop}\label{Prop:bal-hyp_non-existence_Pn-1} Let $X$ be a compact complex manifold with $\mbox{dim}_\C X=n$,\, and $\pi_X:\widetilde{X}\to X$ be the universal covering map of $X$. Suppose that $\pi_X^\star\omega = \partial\alpha+\bar\partial\beta$ for some smooth $(0,1)$-form $\alpha$ and $(1,0)$-form $\beta$ on $\widetilde{X}$ . Then, there is no non-constant holomorphic map $f:\Proj^{1}\longrightarrow X$.

\end{Prop}

  \noindent {\it Proof.}

  Suppose there exists a non-constant holomorphic map $f:\Proj^1\longrightarrow X$ . We will show that $f^\star\omega=0$ on $\Proj^{1}$, contradicting the non-constancy of $f$.

 Let $\tilde{f}:\Proj^{1}\longrightarrow\widetilde{X}$ be a lift of $f$ to $\widetilde{X}$, namely a holomorphic map such that $f=\pi_X\circ\tilde{f}$. From $$f^\star\omega = \tilde{f}^\star(\pi_X^\star\omega),$$ we get by integration: $$ \int\limits_{\Proj^{1}}f^\star\omega = \int\limits_{\Proj^{1}}\tilde{f}^\star(\pi_X^\star\omega) =\int\limits_{\Proj^{1}}\tilde{f}^\star(\partial\alpha+\bar\partial\beta)=\int\limits_{\Proj^{1}}\tilde{f}^\star(d(\alpha+\beta))=\int\limits_{\Proj^{1}}d(\tilde{f}^\star(\alpha+\beta))= 0,$$ where the last identity follows from Stokes's theorem.

Meanwhile, $f^\star\omega\geq 0$ at every point of $\Proj^{1}$. Therefore, $f^\star\omega=0$ on $\Proj^{1}$, a contradiction.  \hfill $\qed$

\vspace{1ex}

An  immediate observation is that, since a {\it SKT hyperbolic} manifold $X$ contains no rational curves, then by Mori’s cone theorem we get $K_X$ is nef.

\subsection{ Gauduchon hyperbolic manifolds}
The second notion that we introduce in this work is the following
\hspace{1ex}
\begin{Def}\label{Def:Gaud-hyperbolic} Let $X$ be a compact complex manifold with $\mbox{dim}_\C X=n$. A Hermitian metric $\omega$ on $X$  is said to be {\bf Gauduchon hyperbolic} if $\omega^{n-1}$ is   $(\widetilde{\partial+\bar\partial})-\mbox{bounded}$ with respect to $\omega$.

   The manifold $X$ is said to be Gauduchon hyperbolic if it carries a {\it Gauduchon hyperbolic metric}.

\end{Def}
Let us first notice the following implication: $$X \hspace{1ex} \mbox{is {\it balanced hyperbolic}} \hspace{1ex} \implies \hspace{1ex} X \hspace{1ex} \mbox{is {\it Gauduchon hyperbolic}}.$$

To see this, besides the obvious fact that every balanced metric is Gauduchon and every $\widetilde{d}$-bounded form is in particular $(\widetilde{\partial+\bar\partial})-\mbox{bounded}$.

 \vspace{1ex}

If $X$ is a compact complex manifold with $\mbox{dim}_\C X=n\geq 2$ and $\omega$ is a Hermitian metric on $X$, for any holomorphic map $f:\C^{n-1}\to X$ that is {\it non-degenerate} at some point $x\in\C^{n-1}$ (in the sense that its differential map $d_xf:\C^{n-1}\longrightarrow T_{f(x)}X$ at $x$ is of maximal rank), we consider the smooth $(1,\,1)$-form $f^\star\omega$ on $\C^{n-1}$. The assumptions made on $f$ imply that the differential map $d_zf$ is of maximal rank for every point $z\in\C^{n-1}\setminus\Sigma$, where $\Sigma\subset\C^{n-1}$ is an analytic subset. Thus, $f^\star\omega$ is $\geq 0$ on $\C^{n-1}$ and is $>0$ on $\C^{n-1}\setminus\Sigma$. Consequently, $f^\star\omega$ can be regarded as a {\it degenerate metric} on $\C^{n-1}$. Its degeneration locus, $\Sigma$, is empty if $f$ is non-degenerate at every point of $\C^{n-1}$, in which case $f^\star\omega$ is a genuine Hermitian metric on $\C^{n-1}$. However, in our case, $\Sigma$ will be non-empty in general, so $f^\star\omega$ will only be a genuine Hermitian metric on $\C^{n-1}\setminus\Sigma$.
Fix an arbitrary integer $n\geq 3$. For any $r>0$, let $B_r:=\{z\in\C^{n-1}\,\mid\,|z|<r\}$ and $S_r:=\{z\in\C^{n-1}\,\mid\,|z|=r\}$ stand for the open ball, resp. the sphere, of radius $r$ centred at $0\in\C^{n-1}$. Moreover, for any $(1,\,1)$-form $\gamma\geq 0$ on a complex manifold and any positive integer $p$, we will use the notation: $$\gamma_p:=\frac{\gamma^p}{p!}.$$
For a holomorphic map $f:\C^{n-1}\to(X,\,\omega)$ in the above setting and for $r>0$, we consider the {\it $(\omega,\,f)$-volume} of the ball $B_r\subset\C^{n-1}$: $$\mbox{Vol}_{\omega,\,f}(B_r):=\int\limits_{B_r}f^\star\omega_{n-1}>0.$$

Meanwhile, for $z\in\C^{n-1}$, let $\tau(z):=|z|^2$ be its squared Euclidean norm. At every point $z\in\C^{n-1}\setminus\Sigma$, we have: \begin{eqnarray}\label{eqn:Hodge-star_area}\frac{d\tau}{|d\tau|_{f^\star\omega}}\wedge\star_{f^\star\omega}\bigg(\frac{d\tau}{|d\tau|_{f^\star\omega}}\bigg) = f^\star\omega_{n-1},\end{eqnarray} where $\star_{f^\star\omega}$ is the Hodge star operator induced by $f^\star\omega$. Thus, the $(2n-3)$-form $$d\sigma_{\omega,\,f}:=\star_{f^\star\omega}\bigg(\frac{d\tau}{|d\tau|_{f^\star\omega}}\bigg)$$ on $\C^{n-1}\setminus\Sigma$ is the area measure induced by $f^\star\omega$ on the spheres of $\C^{n-1}$. This means that its restriction \begin{eqnarray}\label{eqn:area-measure_t_1}d\sigma_{\omega,\,f,\,t}:=\bigg(\star_{f^\star\omega}\bigg(\frac{d\tau}{|d\tau|_{f^\star\omega}}\bigg)\bigg)_{|S_t}\end{eqnarray} is the area measure induced by the degenerate metric $f^\star\omega$ on the sphere $S_t=\{\tau(z)=t^2\}\subset\C^{n-1}$ for every $t>0$. In particular, the area of the sphere $S_r\subset\C^{n-1}$ w.r.t. $d\sigma_{\omega,\,f,\,r}$ is $$A_{\omega,\,f}(S_r)=\int\limits_{S_r}d\sigma_{\omega,\,f,\,r}>0,  \hspace{3ex} r>0.$$
In this case  there is no reason for  the measure on the sphere that "disintegrates" the volume form $f^\star(\omega^{n-1})$, and
the area measure induced by the partially degenerate metric $f^\star(\omega)$  to coincide (Contrary to the case of the equality (b) obtained in \ref{eqn:bal-div-hyperbolic_implication_proof_2_1}, which is valid except when working on $\mathbb{C}$).
In fact, if $\tau (z) = |z|^2$
, the area measure induced by the degenerate metric $f^\star(\omega)$ on $S_t$ is almost by definition
\begin{eqnarray*}\label{eqn:area-measure_t}d\sigma_{\omega,\,f,\,t}:=\bigg(\star_{f^\star\omega}\bigg(\frac{d\tau}{|d\tau|_{f^\star\omega}}\bigg)\bigg)_{|S_t}.\end{eqnarray*}
On the other hand, if $d\mu_{\omega,f,t}$ is the
positive measure that disintegrates $f^\star(\omega^{n-1})$
on the sphere $S_t = \{\tau (z)=t^2\}$, we have
$t = \tau^{\frac{1}{2}}$
, $d\tau = 2t dt$, then we get \begin{eqnarray*}\label{eqn:bal-div-hyperbolic_implication_proof_1_a}\mbox{Vol}_{\omega,\,f}(B_r) = \int\limits_{B_r}f^\star\omega_{n-1} = \int\limits_0^r\bigg(\int\limits_{S_t}d\mu_{\omega,\,f,\,t}\bigg)\,dt = \int_{B_r}d\mu_{\omega,\,f,\,t}\wedge\frac{d\tau}{2t}.\end{eqnarray*}
As a consequence, we get $$(f^\star\omega^{n-1})_{|S_t}=\frac{1}{2}(d\mu_{\omega,\,f,\,t}\wedge d\tau)_{|S_t}.$$

Let us recall the two following definitions introduiced in [MP22a]

\begin{Def}\label{Def:subexp-growth} Let $(X,\,\omega)$ be a {\bf compact} complex Hermitian manifold with $\mbox{dim}_\C X=n\geq 2$ and let $f:\C^{n-1}\to X$ be a holomorphic map that is {\bf non-degenerate} at some point $x\in\C^{n-1}$.
\begin{enumerate}
    \item[1)]  We say that $f$ has {\bf subexponential growth} if the following two conditions are satisfied:

  \vspace{1ex}

  (i)\, there exist constants $C_1>0$ and $r_0> 0$ such that \begin{equation}\label{eqn:comparative-growth}\int\limits_{S_t}|d\tau|_{f^\star\omega}\,d\sigma_{\omega,\,f,\,t}\leq C_1t\,\mbox{Vol}_{\omega,\,f}(B_t),  \hspace{3ex} t>r_0;\end{equation}

\vspace{1ex} (ii)\, for every constant $C>0$, we have: \begin{equation}\label{eqn:subexp-growth}\limsup\limits_{b\to +\infty}\bigg(\frac{b}{C} - \log F(b)\bigg) = +\infty,\end{equation} where $$F(b):=\int\limits_0^b\mbox{Vol}_{\omega,\,f}(B_t)\,dt = \int\limits_0^b\bigg(\int\limits_{B_t}f^\star\omega_{n-1} \bigg)\,dt, \hspace{3ex} b>0.$$
\item[2)] We say that $X$ is {\bf divisorially hyperbolic} if there is no  holomorphic map $f$ that has {\bf subexponential growth}.
\end{enumerate}
\end{Def}

\begin{The}\label{The:Gaud-div-hyperbolic_implication} Every {\bf Gauduchon hyperbolic} compact complex manifold is also {\bf divisorially hyperbolic}.

\end{The}

\noindent {\it Proof.} Let $X$ be a compact complex manifold, with $\mbox{dim}_\C X=n$, equipped with a {\it Gauduchon hyperbolic} metric $\omega$. This means that, if $\pi_X:\widetilde{X}\longrightarrow X$ is the universal cover of $X$, we have $$\pi_X^\star\omega^{n-1} = \partial\Gamma_{n-2,n-1} +\bar\partial\Gamma_{n-1,n-2} \hspace{3ex} \mbox{on}\hspace{1ex} \widetilde{X},$$ where $\Gamma_{n-2,n-1}$ is an $\widetilde\omega$-bounded $C^\infty$ $(n-2,n-1)$-form on $\widetilde{X}$, $\Gamma_{n-1,n-2}$ is an $\widetilde\omega$-bounded $C^\infty$ $(n-1,n-2)$-form on $\widetilde{X}$ and $\widetilde\omega=\pi_X^\star\omega$ is the lift of the metric $\omega$ to $\widetilde{X}$.

Suppose there exists a holomorphic map $f:\C^{n-1}\longrightarrow X$ that is non-degenerate at some point $x\in\C^{n-1}$ and has subexponential growth in the sense of Definition \ref{Def:subexp-growth}. We will prove that $f^\star\omega^{n-1}=0$ on $\C^{n-1}$, in contradiction to the non-degeneracy assumption made on $f$ at $x$.

Since $\C^{n-1}$ is simply connected, there exists a lift $\widetilde{f}$ of $f$ to $\widetilde{X}$, namely a holomorphic map $\tilde{f}:\C^{n-1}\longrightarrow\widetilde{X}$ such that $f=\pi_X\circ\tilde{f}$. In particular, $d_x\tilde{f}$ is injective since $d_xf$ is.

The smooth $(n-1,\,n-1)$-form $f^\star\omega^{n-1}$ is $\geq 0$ on $\C^{n-1}$ and $>0$ on $\C^{n-1}\setminus\Sigma$, where $\Sigma\subset\C^{n-1}$ is the proper analytic subset of all points $z\in\C^{n-1}$ such that $d_zf$ is not of maximal rank. We have: $$f^\star\omega_{n-1} = \tilde{f}^\star(\pi_X^\star\omega^{n-1}) = d(\tilde{f}^\star\Gamma_{n-2,n-1} +\Gamma_{n-1,n-2})  \hspace{3ex} \mbox{on}\hspace{1ex} \C^{n-1}.$$ 

 With respect to the degenerate metric $f^\star\omega$ on $\C^{n-1}$, we have the following

  \begin{Claim}\label{Claim:f-tilde-star_gamma_bounded_0} The $(2n-3)$-form $\tilde{f}^\star(\Gamma_{n-2,n-1} +\Gamma_{n-1,n-2})$ is $(f^\star\omega)$-bounded on $\C^{n-1}$.

  \end{Claim}

  \noindent {\it Proof of Claim.} For any tangent vectors $v_1,\dots , v_{2n-3}$ in $\C^{n-1}$, we have: \begin{eqnarray*}|(\tilde{f}^\star(\Gamma_{n-2,n-1} +\Gamma_{n-1,n-2}))(v_1,\dots , v_{2n-3})|^2 & = & |(\Gamma_{n-2,n-1} +\Gamma_{n-1,n-2})(\tilde{f}_\star v_1,\dots , \tilde{f}_\star v_{2n-3})|^2\\& \stackrel{(a)}{\leq}& C\,|\tilde{f}_\star v_1|^2_{\widetilde\omega}\dots |\tilde{f}_\star v_{2n-3}|^2_{\widetilde\omega} \\
&=&C\,|v_1|^2_{\tilde{f}^\star\widetilde\omega}\dots |v_{2n-3}|^2_{\tilde{f}^\star\widetilde\omega} \stackrel{(b)}{=} C\,|v_1|^2_{f^\star\omega}\dots |v_{2n-3}|^2_{f^\star\omega},\end{eqnarray*} where $C>0$ is a constant independent of the $v_j$'s that exists such that inequality (a) holds thanks to the {\it $\widetilde\omega$-boundedness} of $(\Gamma_{n-2,n-1}$ and $+\Gamma_{n-1,n-2}$ on $\widetilde{X}$, while (b) follows from $\tilde{f}^\star\widetilde\omega = f^\star\omega$.   \hfill $\Box$

  \vspace{2ex}

  \noindent {\it End of Proof of Theorem \ref{The:bal-div-hyperbolic_implication}.} 

\vspace{2ex}

 $\bullet$ On the one hand, we have $d\tau = 2t\,dt$ and \begin{eqnarray}\label{eqn:bal-div-hyperbolic_implication_proof_1_c}\mbox{Vol}_{\omega,\,f}(B_r) = \int\limits_{B_r}f^\star\omega_{n-1} = \int\limits_0^r\bigg(\int\limits_{S_t}d\mu_{\omega,\,f,\,t}\bigg)\,dt = \int_{B_r}d\mu_{\omega,\,f,\,t}\wedge\frac{d\tau}{2t},\end{eqnarray} where $d\mu_{\omega,\,f,\,t}$ is the positive measure on $S_t$ defined by $$\frac{1}{2t}\,d\mu_{\omega,\,f,\,t}\wedge (d\tau)_{|S_t} = (f^\star\omega_{n-1})_{|S_t},  \hspace{3ex} t>0.$$

Comparing with (\ref{eqn:Hodge-star_area}) and (\ref{eqn:area-measure_t_1}), this means that the measures $d\mu_{\omega,\,f,\,t}$ and $d\sigma_{\omega,\,f,\,t}$ on $S_t$ are related by \begin{eqnarray}\label{eqn:sigma-tau_relation}\frac{1}{2t}\,d\mu_{\omega,\,f,\,t} = \frac{1}{|d\tau|_{f^\star\omega}}\,d\sigma_{\omega,\,f,\,t}, \hspace{3ex} t>0.\end{eqnarray}

Now, the H\"older inequality yields: \begin{eqnarray*}\int_{S_t}\frac{1}{|d\tau|_{f^\star\omega}}\,d\sigma_{\omega,\,f,\,t}\geq\frac{A^2_{\omega,\,f}(S_t)}{\int_{S_t}|d\tau|_{f^\star\omega}\,d\sigma_{\omega,\,f,\,t}},\end{eqnarray*} so together with (\ref{eqn:bal-div-hyperbolic_implication_proof_1_c}) and (\ref{eqn:sigma-tau_relation}) this leads to: \begin{eqnarray}\label{eqn:bal-div-hyperbolic_implication_proof_1_b}\nonumber\mbox{Vol}_{\omega,\,f}(B_r) & = & \int\limits_0^r\bigg(\int\limits_{S_t}\frac{1}{2t}\,d\mu_{\omega,\,f,\,t}\bigg)\,d\tau = \int\limits_0^r\bigg(\int\limits_{S_t}\frac{1}{|d\tau|_{f^\star\omega}}\,d\sigma_{\omega,\,f,\,t}\bigg)\,d\tau \\
 & \geq & 2\,\int\limits_0^r\frac{A^2_{\omega,\,f}(S_t)}{\int_{S_t}|d\tau|_{f^\star\omega}\,d\sigma_{\omega,\,f,\,t}}\,t\,dt, \hspace{3ex} r>0.\end{eqnarray}

\vspace{2ex}

 $\bullet$ On the other hand, for every $r>0$, we have: \begin{eqnarray}\label{eqn:bal-div-hyperbolic_implication_proof_2_e}\mbox{Vol}_{\omega,\,f}(B_r) = \int\limits_{B_r} f^\star\omega_{n-1} = \int\limits_{B_r}d(\tilde{f}^\star(\Gamma_{n-2,n-1} +\Gamma_{n-1,n-2}) & = &\int\limits_{S_r}\tilde{f}^\star(\Gamma_{n-2,n-1} +\Gamma_{n-1,n-2})\nonumber \\ &\stackrel{(a)}{\leq} & C\,\int\limits_{S_r}d\sigma_{\omega,\,f} = C\,A_{\omega,\,f}(S_r),\end{eqnarray} where $C>0$ is a constant that exists such that inequality (a) holds thanks to Claim \ref{Claim:f-tilde-star_gamma_bounded_0}.

Putting (\ref{eqn:bal-div-hyperbolic_implication_proof_1_b}) and (\ref{eqn:bal-div-hyperbolic_implication_proof_2_e}) together, we get for every $r>r_0$: \begin{eqnarray}\label{eqn:bal-div-hyperbolic_implication_proof_3}\nonumber\mbox{Vol}_{\omega,\,f}(B_r) & \geq & \frac{2}{C^2}\,\int\limits_0^r\mbox{Vol}_{\omega,\,f}(B_t)\,\frac{t\,\mbox{Vol}_{\omega,\,f}(B_t)}{\int_{S_t}|d\tau|_{f^\star\omega}\,d\sigma_{\omega,\,f,\,t}} \,dt \\
   & \stackrel{(a)}{\geq} & \frac{2}{C_1\,C^2}\int\limits_{r_0}^r\mbox{Vol}_{\omega,\,f}(B_t)\,dt \stackrel{(b)}{:=} C_2\,F(r),\end{eqnarray} where (a) follows from the growth assumption (\ref{eqn:comparative-growth}) and (b) is the definition of a function $F:(r_0,\,+\infty)\longrightarrow(0,\,+\infty)$ with $C_2:=2/(C_1\,C^2)$.

  By taking the derivative of $F$, we get for every $r>r_0$: \begin{eqnarray*}F'(r) = \mbox{Vol}_{\omega,\,f}(B_r) \geq C_2\,F(r),\end{eqnarray*} where the last inequality is (\ref{eqn:bal-div-hyperbolic_implication_proof_3}). This amounts to \begin{eqnarray*}\frac{d}{dt}\bigg(\log F(t)\bigg) \geq C_2, \hspace{3ex} t>r_0.\end{eqnarray*} Integrating this over $t\in[a,\,b]$, with $r_0<a<b$ arbitrary, we get: \begin{eqnarray}\label{eqn:final_ineq_growth-cond}-\log F(a) \geq -\log F(b) + C_2\,(b-a), \hspace{3ex} r_0<a<b.\end{eqnarray}

  Now, fix an arbitrary $a>r_0$ and let $b\to +\infty$. Thanks to the {\it subexponential growth} assumption (\ref{eqn:subexp-growth}) made on $f$, there exists a sequence of reals $b_j\to +\infty$ such that the right-hand side of inequality (\ref{eqn:final_ineq_growth-cond}) for $b=b_j$ tends to $+\infty$ as $j\to +\infty$. This forces $F(a) = 0$ for every $a>r_0$, hence $\mbox{Vol}_{\omega,\,f}(B_r)=0$ for every $r>r_0$. This amounts to $f^\star\omega^{n-1}=0$ on $\C^{n-1}$, in contradiction to the non-degeneracy assumption made on $f$ at a point $x\in\C^{n-1}$.  \hfill $\Box$

\vspace{3ex}
\subsection{ Some properties of SKT hyperbolic manifolds}

It is a classical fact due to Gaffney [Gaf54] that certain basic facts in the Hodge Theory of compact Riemannian manifolds remain valid on {\it complete} such manifolds. The main ingredient in the proof of this fact is the following {\it cut-off trick} of Gaffney's that played a key role in [Gro91, $\S.1$]. It also appears in [Dem97, VIII, Lemma 2.4].

\begin{Lem}([Gaf54])\label{Lem:complete-manifolds_cut-off_functions} Let $(X,\,g)$ be a Riemannian manifold. Then, $(X,\,g)$ is {\bf complete} if and only if there exists an exhaustive sequence $(K_\nu)_{\nu\in\N}$ of {\bf compact} subsets of $X$: $$K_\nu\subset\mathring{K}_{\nu+1} \hspace{3ex} \mbox{for all}\hspace{1ex}\nu\in\N \hspace{3ex} \mbox{and}\hspace{3ex} X = \bigcup\limits_{\nu\in\N}K_\nu,$$ and a sequence $(\psi_\nu)_{\nu\in\N}$ of $C^\infty$ functions $\psi_\nu:X\longrightarrow[0,\,1]$ satisfying, for every $\nu\in\N$, the conditions: \\

\hspace{6ex}  $\psi_\nu = 1$ in a neighbourhood of $K_\nu$,  \hspace{3ex} $\mbox{Supp}\,\psi_\nu\subset\mathring{K}_{\nu+1}$ \hspace{3ex} and

  \vspace{2ex}

  \hspace{6ex}  $||d\psi_\nu||_{L^\infty_g}:=\sup\limits_{x\in X}|(d\psi_\nu)(x)|_g\leq\varepsilon_\nu$,\\

\noindent for some constants $\varepsilon_\nu>0$ such that $\varepsilon_\nu\downarrow 0$ as $\nu$ tends to $+\infty$.    

\end{Lem}  

In particular, the {\it cut-off functions} $\psi_\nu$ are {\it compactly supported}. One can choose $\varepsilon_\nu = 2^{-\nu}$ for each $\nu$ (see e.g. [Dem97, VIII, Lemma 2.4]), but this will play no role here.

\vspace{2ex}

An immediate consequence of Gaffney's cut-off trick is the following classical generalisation of Stokes's Theorem to possibly non-compact, but complete Riemannian manifolds when the forms involved are $L^1$. 

\begin{Lem}([Gro91, Lemma 1.1.A.])\label{Lem:Stokes_non-compact_L1} Let $(X,\,g)$ be a {\bf complete} Riemannian manifold of real dimension $m$. Let $\eta$ be an $L^1_g$-form on $X$ of degree $m-1$ such that $d\eta$ is again $L^1_g$. Then $$\int_Xd\eta = 0.$$

\end{Lem}

\vspace{3ex}
\vspace{2ex}
We will split this section into two parts: one for the K\"ahler case and one for the general case.
\subsection{Case of non-K\"ahlerian manifolds}
\begin{Prop}\label{Prop:no-L1-currents} Let $(X,\,\omega)$ be a {\bf SKT hyperbolic manifold} and let $\pi:\widetilde{X}\longrightarrow X$ be the universal cover of $X$. There exists no non-zero $d$-closed positive $(n-1,\,n-1)$-current $\widetilde{T}\geq 0$ on $\widetilde{X}$ such that $\widetilde{T}$ is $L^1_{\widetilde\omega}$, where $\widetilde\omega:=\pi^\star\omega$ is the lift of $\omega$ to $\widetilde{X}$.

\end{Prop}
\noindent {\it Proof.}
 Let $n=\mbox{dim}_\C X$. The SKT hyperbolic assumption on $X$ means that $\pi^\star\omega = \partial\alpha+\bar\partial\beta$ on $\widetilde{X}$ for some smooth $L^\infty_{\widetilde\omega}$-forms $\alpha$ and $\beta$ of degree $(0,1)$ and $(1,0)$ respectively on $\widetilde{X}$.
 By the closedness assumption on $\widetilde T$ as in the statement, we get 
 \begin{equation}\int\limits_{\widetilde{X}}\widetilde{T}\wedge\pi^\star\omega = \int\limits_{\widetilde{X}}\widetilde{T}\wedge(\partial\alpha+\bar\partial\beta)=\int\limits_{\widetilde{X}}d(\widetilde{T}\wedge(\alpha+\beta)) .\end{equation}
  We would have \begin{equation}\label{eqn:The:no-L1-currents_proof_1}0<\int\limits_{\widetilde{X}}\widetilde{T}\wedge\pi^\star\omega = \int\limits_{\widetilde{X}}d(\widetilde{T}\wedge(\alpha+\beta)) = 0,\end{equation} which is contradictory.

The last identity in (\ref{eqn:The:no-L1-currents_proof_1}) follows from Lemma \ref{Lem:Stokes_non-compact_L1} applied on the complete manifold $(\widetilde{X},\, \widetilde\omega)$ to the $L^1_{\widetilde\omega}$-current $\eta:= \widetilde{T}\wedge(\alpha+\beta)$ of degree $1$ whose differential $d\eta = \widetilde{T}\wedge\pi^\star\omega$ is again $L^1_{\widetilde\omega}$. That $\eta$ is $L^1_{\widetilde\omega}$ follows from $\widetilde T$  $L^1_{\widetilde\omega}$ (by hypothesis) and from $\alpha$ and $\beta$ being both $L^\infty_{\widetilde\omega}$, while $d\eta$ being $L^1_{\widetilde\omega}$ follows from $\widetilde{T}$ being $L^1_{\widetilde\omega}$ and from $\pi^\star\omega$ being $L^\infty_{\widetilde\omega}$ (as a lift of the smooth, hence bounded, form $\omega$ on the {\it compact} manifold $X$).

 \hfill $\Box$
 
 \vspace{2ex}
 Similarly, we can replace the SKT hyperbolic manifold with a Gauduchon hyperbolic manifold, and using the same approach, we can demonstrate that there does not exist a non-zero $d$-closed positive $(1,1)$-current $\widetilde{T}\geq 0$ on $\widetilde{X}$ such that $\widetilde{T}$ is $L^1_{\widetilde\omega}$.
\begin{The}\label{The:bal-complete-exact_no-harmonic-forms_degree1} Let $X$ be a compact complex {\bf SKT hyperbolic} manifold with $\mbox{dim}_\C X=n$. Let $\pi:\widetilde{X}\longrightarrow X$ be the universal cover of $X$ and $\widetilde\omega:=\pi^\star\omega$ the lift to $\widetilde{X}$ of a SKT hyperbolic metric $\omega$ on $X$. Fix a primitive $L_{\widetilde{\omega}}^2$-form $\phi$ on $\widetilde{X}$  of bidegree $(p,q)$ with $p+q=n-1$ such that $$ \partial\phi=0, \quad\quad\bar\partial\phi=0.$$ Then $\phi=0$.
\end{The}
\noindent {\it Proof.}
Let $n=\mbox{dim}_\C X$. The SKT hyperbolic assumption on $X$ means that $\pi^\star\omega = \partial\alpha+\bar\partial\beta$ on $\widetilde{X}$ for some smooth $L^\infty_{\widetilde\omega}$-forms $\alpha$ and $\beta$ of degree $(0,1)$ and $(1,0)$ respectively on $\widetilde{X}$.
Let $\phi$ be as in the statement of the theorem. Then, by \ref{complete_case} of Lemma \ref{Lem:1-forms_Delta-harm} we get 
\begin{eqnarray}\label{0_with_partial}
    0=\big<\big<\partial^\star(\widetilde{\omega}\wedge\phi),\psi_\nu\alpha\wedge\phi\big>\big>=\big<\big<\widetilde{\omega}\wedge\phi,\partial(\psi_\nu\alpha\wedge\phi)\big>\big>&=&\big<\big<\widetilde{\omega}\wedge\phi,\partial\psi_\nu\wedge\alpha\wedge\phi\big>\big>\nonumber\\
    &+& \big<\big<\widetilde{\omega}\wedge\phi,\psi_\nu\wedge\partial\alpha\wedge\phi\big>\big> ,
\end{eqnarray}
and  
\begin{eqnarray}\label{0_with_bar_partail}
    0=\big<\big<\bar\partial^\star(\widetilde{\omega}\wedge\phi),\psi_\nu\beta\wedge\phi\big>\big>=\big<\big<\widetilde{\omega}\wedge\phi,\bar\partial(\psi_\nu\beta\wedge\phi)\big>\big>&=&\big<\big<\widetilde{\omega}\wedge\phi,\bar\partial\psi_\nu\wedge\beta\wedge\phi\big>\big>\nonumber\\
    &+& \big<\big<\widetilde{\omega}\wedge\phi,\psi_\nu\wedge\bar\partial\beta\wedge\phi\big>\big> .
\end{eqnarray}
Putting \ref{0_with_partial} and \ref{0_with_bar_partail} together, we get 
\begin{eqnarray}\label{both_0}
    \big<\big<\widetilde{\omega}\wedge\phi,d(\psi_\nu)\wedge(\alpha+\beta)\wedge\phi\big>\big>+\big<\big<\widetilde{\omega}\wedge\phi,\psi_\nu\wedge(\partial\alpha+\bar\partial\beta)\wedge\phi\big>\big>=0.
\end{eqnarray}
Therefore, 
\begin{eqnarray}\label{premier_terme_nul}
    \bigg|\int\limits_{\widetilde{X}}<\big<\widetilde{\omega}\wedge\phi,d(\psi_\nu)\wedge(\alpha+\beta)\wedge\phi\big>dV\bigg|&\leq& \int\limits_{\widetilde{X}}||\widetilde{\omega}||_{L^\infty}||d\psi_\nu||_{L^\infty}||\alpha+\beta||_{L^\infty}||\phi||_{L^2}dV\nonumber\\&\leq&\varepsilon_\nu||\widetilde{\omega}||_{L^\infty}||\alpha+\beta||_{L^\infty}||\phi||_{L^2} \downarrow 0 \mbox{ as } \nu \to +\infty.
\end{eqnarray}
Hence, combining \ref{premier_terme_nul} and  \ref{both_0} we get,
\begin{eqnarray}
\lim_{\nu\to+\infty}\big<\big<\widetilde{\omega}\wedge\phi,\psi_\nu\wedge(\partial\alpha+\bar\partial\beta)\wedge\phi\big>\big>=\big<\big<\widetilde{\omega}\wedge\phi,\widetilde\omega\wedge\phi\big>\big>=0
\end{eqnarray}
On the other hand, the pointwise map $\widetilde\omega\wedge\cdot : \Lambda^{n-1}T^\star \widetilde{X}\longrightarrow\Lambda^{n+1}T^\star \widetilde{X}$ is bijective, so we get $\phi=0$ from $\widetilde\omega\wedge \phi=0$.\\
The proof is complete

 \hfill $\Box$
 \begin{Cor}
  Let
$\phi$ be a $(n-1, 0)$-form (respectively a $(0,n-1)$-form) on a connected complete manifold  $(\widetilde{X},\widetilde{\omega})$ such that
$$\phi\in L^2(\widetilde{X}), \quad\quad \partial\phi= 0, \quad\quad \bar\partial\phi = 0.$$
If $\widetilde{\omega}=\partial\alpha+\bar\partial\beta$ where $\alpha$ and $\beta$ are bounded $1$-forms on $\widetilde{X}$, then $$\phi=0.$$  
 \end{Cor}

 Similarly, since all $1$-forms are primitive, we get the following
 \begin{Prop}
     Let
$\phi$ be a $(1, 0)$-form (respectively a $(0,1)$-form) on a connected complete manifold  $(\widetilde{X},\widetilde{\omega})$ such that
$$\phi\in L^2(\widetilde{X}), \quad\quad \partial\phi= 0, \quad\quad \bar\partial\phi = 0.$$
If $\widetilde{\omega}^{n-1}=\partial\alpha+\bar\partial\beta$ where $\alpha$ and $\beta$ are bounded $(2n-3)$-forms on $\widetilde{X}$, then $$\phi=0.$$ 
 \end{Prop}
\subsection{Case of K\"ahler manifolds}
Let $(X,\omega)$ be a  K\"ahler manifold of dimension $2n$. The operator $L^k: \Lambda^pT^\star X\to \Lambda^{2k+p}T^\star X$ defined by $L^k(\phi)=\omega^k\wedge\phi$ for all $p$-forms $ \phi\in\Lambda^pT^\star X$
 commutes with $d$ and $\Delta$.
We first recall the following result
\begin{The}(Lefschetz.) The map $L^k$ is injective on harmonic forms for $2p + 2k \leq 2n = dim_\R(X)$ and surjective for $2p + 2k\geq 
 2n$.
\end{The}
This result show that $L^k$
 is bijective for $2p + 2k =
2n$. Moreover, the injectivity (in particular the bijectivity) shows that $L$ is a {\bf quasi-isometry}, i.e., $$ C^{-1}||\phi||_{L^2}\leq||L^k\phi||_{L^2}\leq C||\phi||_{L^2} \quad\quad\mbox{for all } \phi\in\Lambda^pT^\star X \mbox{ and some constant }C>0.$$
For more details see Theorem $1.2.A$. in [Gro91].
\begin{The}
Let $(X, \omega)$ be a complete K\"ahler manifold of dimension $2n $ and $\omega=\partial\alpha+\bar\partial\beta$ where $\alpha$ and $\beta$ are respectively a bounded $(0,1)$ and $(1,0)$ forms   on $X$. Then every $L^2$-form $\psi$ on $X$ of degree $p\neq n$ satisfies the inequality
$$\big<\psi,\Delta\psi\big>\geq\lambda_0^2\big<\psi,\psi\big>,$$
where $\lambda_0$
 is a strictly positive constant which depends only on $n = dim_\C X,\, \alpha$ and $\beta$.
 \end{The}
\noindent {\it Proof.}
 For a given $p < n$, let $k=n-p$, we get $2p+2k=2n$, so by the Lefschetz theorem $L^k$ is a bijective quasi-isometry and so every $L^2$
-form $\psi$ of degree $2k+p$ is the product $\psi=L^k\phi=\omega^k\wedge\phi$, where $\phi=L^{-k}\psi$ and $C^{-1}||\phi||_{L^2}\leq||L^k\phi||_{L^2}\leq C||\phi||_{L^2} $. Since $L^k$ commutes with $\Delta$, we also have $$ \big<\big<\Delta\phi,\phi\big>\big>\leq C_1 \big<\big<\Delta\psi,\psi\big>\big>$$
Meanwhile, we have: $\psi= \omega^k\wedge\phi = (\partial\alpha+\bar\partial\beta)\wedge\omega^{k-1}\wedge\phi$. In other words, \begin{eqnarray}\label{eqn:bal-complete-exact_Laplacian_l-bound_3}\psi= d\theta -\psi', \hspace{2ex} \mbox{where}\hspace{1ex} \theta:=(\alpha+\beta)\wedge\omega^{k-1}\wedge\phi \hspace{2ex}\mbox{and}\hspace{2ex} \psi':=(\alpha+\beta)\wedge\omega^{k-1}\wedge d\phi+(\bar\partial\alpha+\partial\beta)\wedge\omega^{k-1}\wedge\phi.\end{eqnarray}

To estimate $\theta$, we write: \begin{eqnarray}\label{eqn:bal-complete-exact_Laplacian_l-bound_4}||\theta||\leq||\alpha+\beta||_{L^\infty_\omega}\,||\phi||_{L^2} \leq C||\alpha+\beta||_{L^\infty_\omega}\,||\psi||_{L^2},\end{eqnarray}.

 To estimate the left part of $\psi'$, since
$||d\phi||_{L^2}^2\leq\big<\big<\Delta\phi,\phi\big>\big>\leq C_1 \big<\big<\Delta\psi,\psi\big>\big>, $

we have: \begin{eqnarray}\label{eqn:bal-complete-exact_Laplacian_l-bound_5}||(\alpha+\beta)\wedge\omega^{k-1}\wedge d\phi||_{L^2}\leq||\alpha+\beta||_{L^\infty_\omega}\,||d \phi||_{L^2} \leq ||\alpha+\beta||_{L^\infty_\omega}\,\langle\langle\Delta\psi,\,\psi\rangle\rangle^{\frac{1}{2}},
\end{eqnarray}  
 
 To find an upper bound for $||\psi||_{L^2}$, we write: \begin{eqnarray}\label{eqn:bal-complete-exact_Laplacian_l-bound_6}||\psi||^2 = \langle\langle\psi,\,d\theta -\psi'\rangle\rangle\leq|\langle\langle\psi,\,d\theta\rangle\rangle| +|\langle\langle\psi,\,\psi'\rangle\rangle|,\end{eqnarray} where (\ref{eqn:bal-complete-exact_Laplacian_l-bound_3}) was used to get the first equality.

For the first term on the r.h.s. of (\ref{eqn:bal-complete-exact_Laplacian_l-bound_6}), we get: \begin{eqnarray}\label{eqn:bal-complete-exact_Laplacian_l-bound_7}|\langle\langle\psi,\,d\theta\rangle\rangle|=|\langle\langle d^\star\psi,\,\theta\rangle\rangle|\leq||d^\star\psi||_{L^2}\,||\theta||_{L^2}\leq C\langle\langle\Delta\psi,\,\psi\rangle\rangle^{\frac{1}{2}}\,||\alpha+\beta||_{L^\infty_\omega}\,||\psi||_{L^2},\end{eqnarray} where (\ref{eqn:bal-complete-exact_Laplacian_l-bound_4}) was used to get the last inequality.

For the second term on the r.h.s. of (\ref{eqn:bal-complete-exact_Laplacian_l-bound_6}),  we get: \begin{eqnarray}\label{eqn:bal-complete-exact_Laplacian_l-bound_8}|\langle\langle\psi,\,\psi'\rangle\rangle| &\leq &|\langle\langle\psi,\,(\alpha+\beta)\wedge\omega^{k-1}\wedge d\phi\rangle\rangle| +|\langle\langle\psi,\,\bar\partial\alpha\wedge\omega^{k-1}\wedge\phi\rangle\rangle|+|\langle\langle\psi,\,\partial\beta\wedge\omega^{k-1}\wedge\phi\rangle\rangle| \nonumber\\ &\leq& ||\psi||_{L^2}\,||\alpha+\beta||_{L^\infty_\omega}\,\langle\langle\Delta\psi,\,\psi\rangle\rangle^{\frac{1}{2}}+|\langle\langle\psi,\,\bar\partial\alpha\wedge\omega^{k-1}\wedge\phi\rangle\rangle|+|\langle\langle\psi,\,\partial\beta\wedge\omega^{k-1}\wedge\phi\rangle\rangle|,\end{eqnarray}  where (\ref{eqn:bal-complete-exact_Laplacian_l-bound_5}) was used to get the last inequality.\\
To estimate the last two terms on the right in \ref{eqn:bal-complete-exact_Laplacian_l-bound_8}, we write:
\begin{eqnarray}\label{eqn:bal-complete-exact_Laplacian_l-bound_9}
    |\langle\langle\psi,\,\bar\partial\alpha\wedge\omega^{k-1}\wedge\phi\rangle\rangle|&\leq&|\langle\langle\psi,\,\bar\partial(\alpha\wedge\omega^{k-1}\wedge\phi)\rangle\rangle|+|\langle\langle\psi,\,\alpha\wedge\omega^{k-1}\wedge\bar\partial\phi\rangle\rangle|\nonumber\\
    &=& |\langle\langle\bar\partial^\star\psi,\,\alpha\wedge\omega^{k-1}\wedge\phi\rangle\rangle|+|\langle\langle\psi,\,\alpha\wedge\omega^{k-1}\wedge\bar\partial\phi\rangle\rangle|\nonumber \\
    &\leq& C||\bar\partial^\star\psi||_{L^2}||\alpha||_{L^\infty_\omega}||\psi||_{L^2}+C||\bar\partial\phi||_{L^2}||\alpha||_{L^\infty_\omega}||\psi||_{L^2}\nonumber\\ &\leq& 2C.C_2||\alpha||_{L^\infty_\omega}||\psi||_{L^2}\big<\big<\Delta\psi,\psi\big>\big>^{\frac{1}{2}}
\end{eqnarray}
since
$||\bar\partial^\star\psi||_{L^2}^2\leq\big<\big<\Delta\phi,\phi\big>\big>\leq C_2 \big<\big<\Delta\psi,\psi\big>\big>, $ and $||\bar\partial\phi||_{L^2}^2\leq\big<\big<\Delta\phi,\phi\big>\big>\leq C_2 \big<\big<\Delta\psi,\psi\big>\big>$.
\vspace{3ex}
Using the same technique, we obtain the following estimate 
\begin{eqnarray}\label{eqn:bal-complete-exact_Laplacian_l-bound_10}
|\langle\langle\psi,\,\partial\beta\wedge\omega^{k-1}\wedge\phi\rangle\rangle|\leq 2C.C_2||\beta||_{L^\infty_\omega}||\psi||_{L^2}\big<\big<\Delta\psi,\psi\big>\big>^{\frac{1}{2}}
\end{eqnarray}

 Adding up (\ref{eqn:bal-complete-exact_Laplacian_l-bound_7}), (\ref{eqn:bal-complete-exact_Laplacian_l-bound_8}), (\ref{eqn:bal-complete-exact_Laplacian_l-bound_9})
  and (\ref{eqn:bal-complete-exact_Laplacian_l-bound_10}) and  using (\ref{eqn:bal-complete-exact_Laplacian_l-bound_6}), we get \begin{eqnarray*}||\psi||^2\leq (\lambda_0^2)^{-1}\,\langle\langle\Delta\psi,\,\psi\rangle\rangle. \end{eqnarray*}
 for the forms $\psi$ of degree $2k+p > n$.
 \vspace{3ex}
 
 The case $2k+p < n$ follows by the
Poincare duality as the operator $\star:\,\Lambda^pT^\star X\to\Lambda^{n-p}T^\star X$ commutes with $\Delta$ and is isometric for the $L_2
$-norms.
 
 The proof is complete.
 
 \hfill $\Box$
 
 As a simple consequence we get the following {\bf Lefschetz vanishing theorem}.
\begin{Cor}  Let $(\widetilde{X},\widetilde{\omega})$ be a connected complete K\"ahler manifold. If $\widetilde{\omega}=\partial\alpha+\bar\partial\beta$ where $\alpha$ and $\beta$ are bounded $1$-forms on $\widetilde{X}$, then ${\cal H}^p_{\Delta_{\widetilde\omega}}(\widetilde{X},\,\C) =  0$, unless $p = n$.

\end{Cor}
\noindent {\bf References.} \\

\vspace{1ex}



\noindent [Gau77a]\, P. Gauduchon --- {\it Le th\'eor\`eme de l'excentricit\'e nulle} --- C. R. Acad. Sci. Paris, S\'er. A, {\bf 285} (1977), 387-390.

\vspace{1ex}

\noindent [Gro91]\, M. Gromov --- {\it K\"ahler Hyperbolicity and $L^2$ Hodge Theory} --- J. Diff. Geom. {\bf 33} (1991), 263-292.

\vspace{1ex}

\noindent [Lam99]\, A. Lamari --- {\it Courants k\"ahl\'eriens et surfaces compactes} --- Ann. Inst. Fourier {\bf 49}, no. 1 (1999), 263-285.

\vspace{1ex}

\noindent [BDPP04]\, S. Boucksom,  J.P. Demailly , M. Paun, , Peternell, T --- {\it The pseudo-effective cone of a compact K\"ahler manifold and varieties of negative Kodaira dimension.} --- arXiv preprint math/0405285, 2004.

\vspace{1ex}

\noindent [LT93]\, P. Lu, G. Tian -- {\it The Complex Structures on Connected Sums of $S^3\times S^3$} --- in Manifolds and Geometry (Pisa, 1993), Sympos. Math., XXXVI, Cambridge Univ. Press, Cambridge, 1996, p. 284–293.

\vspace{1ex}

\noindent [Mic83]\, M. L. Michelsohn --- {\it On the Existence of Special Metrics in Complex Geometry} --- Acta Math. {\bf 143} (1983) 261-295.

\vspace{1ex}
\noindent [MP22a]\, S. Marouani, D. Popovici --- {\it Balanced Hyperbolic and Divisorially Hyperbolic Compact Complex Manifolds } ---arXiv e-print CV 2107.08972v2, to appear in Mathematical Research Letters.

\vspace{1ex}

\noindent [MP22b]\, S. Marouani, D. Popovici --- {\it Some properties of balanced hyperbolic
compact complex manifolds } --- Internat. J. Math. 33 (2022), no. 3, Paper No.
2250019, 39. MR 4390652.

\vspace{1ex}
\noindent [Pop13]\, D. Popovici --- {\it Deformation Limits of Projective Manifolds: Hodge Numbers and Strongly Gauduchon Metrics} --- Invent. Math. {\bf 194} (2013), 515-534.

\vspace{1ex}

\noindent [Pop15]\, D. Popovici --- {\it Aeppli Cohomology Classes Associated with Gauduchon Metrics on Compact Complex Manifolds} --- Bull. Soc. Math. France {\bf 143}, no. 4 (2015), p. 763-800.
\vspace{1ex}

\noindent [PT]\, R. Piovani, A. Tomassini,  --- {\it Aeppli cohomology and Gauduchon metrics} --- Complex Analysis and Operator Theory, 2020, vol. 14, no 1, p. 22.

\vspace{1ex}

\noindent [PU18]\, D. Popovici, L. Ugarte --- {\it Compact Complex Manifolds with Small Gauduchon Cone} --- Proc. LMS {\bf 116}, no. 5 (2018) doi:10.1112/plms.12110.

\vspace{1ex}

\noindent [Sch07]\, M. Schweitzer --- {\it Autour de la cohomologie de Bott-Chern} --- arXiv e-print math.AG/0709.3528v1.
\vspace{1ex}
\noindent [YZZ22]\,S.T. Yau, Q. Zhao, F. Zheng, On Strominger K\"ahler-like manifolds with degenerate torsion,
preprint arXiv:1908.05322 [math.DG].
\vspace{1ex}

\noindent [ST10]\, J. Streets, G. Tian --- {\it A Parabolic Flow of Pluriclosed Metrics} --- Int. Math. Res. Notices, {\bf 16} (2010), 3101-3133.

\vspace{1ex}

\noindent [Voi02]\, C. Voisin --- {\it Hodge Theory and Complex Algebraic Geometry. I.} --- Cambridge Studies in Advanced Mathematics, 76, Cambridge University Press, Cambridge, 2002.



\vspace{1ex}

\noindent [YZZ19]\, S.-T. Yau, Q. Zhao, F. Zheng --- {\it On Strominger K\"ahler-like Manifolds with Degenerate Torsion} --- arXiv e-print DG 1908.05322v2

\vspace{3ex}

\noindent Institut de Math\'ematiques de Toulouse, Universit\'e Paul Sabatier,

\noindent 118 route de Narbonne, 31062 Toulouse, France

\noindent Email 1: 	almarouanisamir@gmail.com
\noindent Email 2: 	samir.marouani@math.univ-toulouse.fr 
\end{document}